%% file: main.tex
\documentclass[11pt,a4paper]{article}
\pdfoutput=1
\input{_imports}
\input{_commands}
\bibliography{bibliography.bib} 
\author{Tim Höpfner\footnote{This work is part of the author's doctorial thesis at the University of Göttingen.}}
\title{Two-Parameter Novikov-Shubin Invariants\\for Fibre Bundles}
\date{\today}
	
\begin{document} 
\maketitle
\begin{abstract} 
	\input{Abstract} 
\end{abstract}
\tableofcontents 
\section{Introduction}
	\input{Introduction}
\section{Two-Parameter Novikov-Shubin Numbers}
	\subsection{Scaling of Fibre Bundles}
		\input{Definition}
	\subsection{Near Cohomological Approach} 
		\input{NearCohom}
\section{Invariance Properties} 
	\input{Invariance}
	\subsection{Metric Invariance for Fixed Connection} 
		\input{FixedConnMetricInvar}
	\subsection{Fibre Homotopy Invariance} 
		\input{FibreHomotopyInvar}
	\subsection{(Partial) Metric Invariance} 
		\input{MetricInvar}
\section{Example: The Heisenberg Group} 
	\input{Example} 
\printbibliography 
\end{document}

%% file: _imports.tex
\usepackage[utf8]{inputenc}
\usepackage{amsmath}
\usepackage{amsfonts}
\usepackage{amssymb}
\usepackage{amsthm}

\usepackage{geometry}
\usepackage{mathtools}

\usepackage{hyperref}
\hypersetup{colorlinks=true,allcolors=.,urlcolor=blue}

\usepackage{tikz-cd}
\usetikzlibrary{decorations.pathmorphing, arrows, decorations.markings}

\usepackage{graphicx}
\usepackage{float}

\usepackage{units} 
\usepackage{fancyvrb} 

\usepackage{pgfplots}
\pgfplotsset{compat=1.15}
\usepackage{mathrsfs}

\usepackage[backend=biber, style=alphabetic]{biblatex}

\usepackage{ifthen}

%% file: _commands.tex
\theoremstyle{plain}
\newtheorem{thm}{Theorem}
\numberwithin{thm}{section}
\newtheorem*{thm*}{Theorem}
\newtheorem{prp}[thm]{Proposition}
\newtheorem{lma}[thm]{Lemma}
\newtheorem{crl}[thm]{Corollary}
\newtheorem{kor}[thm]{Korollar}

\theoremstyle{definition}
\newtheorem{dfn}[thm]{Definition}
\newtheorem{exm}[thm]{Example}
\newtheorem{bsp}[thm]{Beispiel}
\newtheorem{exer}[thm]{exercise}
\newtheorem{aufg}[thm]{Aufgabe}
\newtheorem{rmk}[thm]{Remark}
\newtheorem{nta}[thm]{Notation}
\newtheorem{asp}[thm]{Assumption}

\newcommand{\theorem}[2][]{
\ifthenelse{ \equal{#1}{} }
      {\begin{thm}{#2}\end{thm}}
      {\begin{thm}[#1]{#2}\end{thm}}
}

\newcommand{\nntheorem}[2][]{
\ifthenelse{ \equal{#1}{} }
      {\begin{thm*}{#2}\end{thm*}}
      {\begin{thm*}[#1]{#2}\end{thm*}}
}

\newcommand{\proposition}[2][]{
\ifthenelse{ \equal{#1}{} }
      {\begin{prp}{#2}\end{prp}}
      {\begin{prp}[#1]{#2}\end{prp}}
}
\newcommand{\lemma}[2][]{
\ifthenelse{ \equal{#1}{} }
      {\begin{lma}{#2}\end{lma}}
      {\begin{lma}[#1]{#2}\end{lma}}
}
\newcommand{\corollary}[2][]{
\ifthenelse{ \equal{#1}{} }
      {\begin{crl}{#2}\end{crl}}
      {\begin{crl}[#1]{#2}\end{crl}}
}
\newcommand{\korollar}[2][]{
\ifthenelse{ \equal{#1}{} }
      {\begin{kor}{#2}\end{kor}}
      {\begin{kor}[#1]{#2}\end{kor}}
}

\newcommand{\definition}[2][]{
\ifthenelse{ \equal{#1}{} }
      {\begin{dfn}{#2}\end{dfn}}
      {\begin{dfn}[#1]{#2}\end{dfn}}
}
\newcommand{\example}[2][]{
\ifthenelse{ \equal{#1}{} }
      {\begin{exm}{#2}\end{exm}}
      {\begin{exm}[#1]{#2}\end{exm}}
}
\newcommand{\beispiel}[2][]{
\ifthenelse{ \equal{#1}{} }
      {\begin{bsp}{#2}\end{bsp}}
      {\begin{bsp}[#1]{#2}\end{bsp}}
}
\newcommand{\exercise}[2][]{
\ifthenelse{ \equal{#1}{} }
      {\begin{exer}{#2}\end{exer}}
      {\begin{exer}[#1]{#2}\end{exer}}
}
\newcommand{\aufgabe}[2][]{
\ifthenelse{ \equal{#1}{} }
      {\begin{aufg}{#2}\end{aufg}}
      {\begin{aufg}[#1]{#2}\end{aufg}}
}
\newcommand{\remark}[2][]{
\ifthenelse{ \equal{#1}{} }
      {\begin{rmk}{#2}\end{rmk}}
      {\begin{rmk}[#1]{#2}\end{rmk}}
}
\newcommand{\notation}[2][]{
\ifthenelse{ \equal{#1}{} }
      {\begin{nta}{#2}\end{nta}}
      {\begin{nta}[#1]{#2}\end{nta}}
}
\newcommand{\assumption}[2][]{
\ifthenelse{ \equal{#1}{} }
      {\begin{asp}{#2}\end{asp}}
      {\begin{asp}[#1]{#2}\end{asp}}
}

\renewcommand{\cal}[1]{\mathcal{#1}}
\renewcommand{\frak}[1]{\mathfrak{#1}}

\newcommand{\N}{\mathbb{N}}
\newcommand{\Z}{\mathbb{Z}}
\newcommand{\R}{\mathbb{R}}

\newcommand{\acts}{\curvearrowright}

\newcommand{\set}[2]{\left\lbrace #1 \:\middle|\: #2 \right\rbrace}

\newcommand{\up}{\mathrm{up}}

\DeclareMathOperator{\im}{im}

\DeclareMathOperator{\tr}{tr}

\DeclareMathOperator{\id}{id}

\DeclareMathOperator{\rk}{rk}

\newcommand{\half}{\nicefrac{1}{2}}
\newcommand{\leftquotient}[2]{{\left.\raisebox{-.15em}{$#1$}\middle\backslash\raisebox{.05em}{$#2$}\right.}}

%% file: Abstract.tex
In this paper we construct a two-parameter version of spectral density functions and Novikov-Shubin invariants on fibre bundles. 
The aim of this approach is to gain a better understanding of how the near-zero spectrum of the Hodge Laplace operators on the fibre and the base of a fibre bundle contribute separately to the near-zero spectrum of the Laplace operators of the total space.
We show that this two-parameter generalisation of the classical spectral density function still satisfies several invariance properties. 
As an example, we compute it explicitly for the three-dimensional Heisenberg group.

%% file: Introduction.tex
Given a product space $M=M_1\times M_2$, we have a good understanding of the near-zero spectrum of the Hodge Laplace operators on $M$ in terms of the near-zero spectra on $M_1$ and $M_2$.
Indeed, in this case we get results based on the Künneth formula, see e.g., the book of W.~Lück~\cite[§2.1.3]{LckL2}.
Similar approaches for non-trivial fibre bundles, for example using the Serre spectral sequence, do not seem to work out as nicely. 

To understand this problem better in the case of fibre bundles, we define two-parameter versions $\cal G_k\colon \R_+\times \R_+\to [0,\infty]$ of spectral density functions and the Novikov-Shubin invariants.
The aim of this generalisation of the classical spectral density functions $\cal F_k$ of the total space is to detect the individual contributions from the base and from the fibre and to obtain finer invariants for (non-trivial) bundles. 
We prove several invariance properties of these generalisations.
First, we show that for a fixed connection the numbers are invariant under change of compatible metrics:
\nntheorem[\ref{Thm_2NSI_Metric_fixed_conn_invar}]{
    Let $G\acts (M\to B,\nabla,g)$ be a fibre bundle with fixed connection $\nabla$ and compatible free proper cocompact group action by a group $G$. Then the dilatational equivalence class of the spectral density function underlying the two-parameter Novikov-Shubin numbers
    $$ \cal G_k(M\to B, \nabla) = \cal G_k(M\to B, \nabla, g)$$
    does not depend on the choice of $G$-invariant $\nabla$-compatible Riemannian metric $g$. 
}

Then, we prove that it is further invariant under certain compatible fibre homotopy equivalences:
\nntheorem[\ref{Thm_2NSI_fibre_homot_invar}]{
If there is a $G$-equivariant fibre homotopy equivalence between suitable bundles $M\to B$ and $M'\to B$ such that $\nabla = f^*\nabla'$, then their spectral density functions are dilatationally equivalent,
$$\cal G_k(M'\to B, \nabla') \sim \cal G_k(M\to B, f^*\nabla').$$
}

Lastly, we show that the two-parameter Novikov-Shubin numbers are invariant under change of connection as long as the fibre is shrunk at least as fast as the base:

\nntheorem[\ref{Thm_2NSI_half_metric_indep}]{
Let $G$ be a group and $M\to B$ be equipped with two pairs of connection and compatible Riemannian metric such that $G\acts (M\to B,\nabla,g)$ and $G\acts (M\to B,\nabla',g')$ are Riemannian fibre bundles with connection and compatible free proper cocompact $G$-action. 
Then the two-parameter spectral density functions restricted to the subspace $\{\nu\leq\mu\}$ are dilatationally equivalent,
$$ \cal G_k(M,\nabla,g)|_{\{\nu\leq\mu\}}\sim \cal G_k(M,\nabla',g')|_{\{\nu\leq\mu\}}. $$ 
}

Then, as an example, we compute these numbers explicitly in the example of the three-dimensional Heisenberg group. 

This fits into the recent study of fibre bundles by means of invariants, such as work based on J.-M.~Bismut and J.~Cheeger's study~\cite{BismutCheeger} of higher torsion invariants on fibre bundles using adiabatic limits, J.-M.~Bismut's study~\cite{Bismut} of an Atiyah-Singer theorem for families of Dirac operators and
characteristic classes of fibre bundles such as the Morita-Miller-Mumford classes named after D. Mumford~\cite{Mumford}, E.~Y.~Miller~\cite{Miller} and S.~Morita~\cite{Morita}.
This new definition is different, but similar in spirit, to the study of adiabatic limits of fibre bundles.

%% file: Definition.tex
Let $(M,g)$ be a non-compact Riemannian manifold with a cocompact free proper group action $G\acts M$ acting by isometries.
The spectral density function $\cal F(d)$ of $d$ is defined in terms of the upper Laplacian $\Delta^k_\up = d^*d\acts  L^2\Omega^k(M,g)$ using the von Neumann trace of the group von Neumann algebra $\cal NG$ of $G$ by
$$ \cal F_k(M,g)(\lambda) = \tr_{\cal NG} \chi_{[0,\lambda^2]}(\Delta^k_\up(M,g))
= \dim_{\cal NG} \im \chi_{[0,\lambda^2]}(\Delta^k_\up(M,g)).$$

An interesting observation is that it can instead be defined by rescaling the manifold and looking at a fixed interval of the spectrum. 
For $g_\lambda= \lambda^2 g$ the Laplace operators for the different metrics satisfy $\Delta^k_\up(M, g_\lambda) = \lambda^{-2} \Delta^k_\up(M,g)$.
Therefore, we obtain 
$$ \cal F_k(M,g)(\lambda) =  \tr_{\cal NG} \chi_{[0,\lambda^2]}(\Delta^k_\up(M,g)) = \tr_{\cal NG} \chi_{[0,1]}(\Delta_\up^k(M, g_\lambda)) = \cal F_k(M, g_\lambda)(1).$$
If $M$ is the total space of a fibre bundle and the fibre bundle structure is compatible with the $G$-action, then we can scale $M$ with different speed in fibre and base directions. 
This way we can define a refined version of the Novikov-Shubin invariants.

More precisely, let 
$M \xrightarrow{\pi} B $
be a fibre bundle with fibres $\{F_b = \pi^{-1}(b)\}_{b\in B}$, where $(M,g)$ is a Riemannian manifold with  Riemannian metric $g$. 
(Without loss of generality, we can assume the base, the fibres and thus also the total space to be connected.)
At every point $x\in M$ with $\pi(x)=b$, we have the subspace
$$ T_x F_b = \ker D_x\pi \subset T_x M, $$
giving rise to the vertical subbundle $VM\subset TM$ of the tangent bundle by
$$ VM = T_\bullet F_\bullet = \ker(\pi_\ast) \subset TM. $$
Choosing a connection $\nabla$ compatible with $g$ on the fibre bundle is equivalent to specifying an orthogonal complement $HM$ of $VM$ in $TM$, so that the tangent space $TM$ decomposes as
$$ TM \cong_\nabla VM \perp_g HM. $$
The bundle $HM$ is called the horizontal subbundle of $TM$.
The Riemannian metric decomposes fibrewise into 
a vertical and a horizontal contribution,
$$ g_x = g_{x,V} + g_{x,H}, $$
where $g_{x,V}$ is supported in $V_xM\otimes V_xM$ and $g_{x,H}$ is supported in $H_xM\otimes H_xM$. 

In the following, we denote the situation described here by the triple $(M\to B, \nabla, g)$ and call such a triple a Riemannian fibre bundle with connection.

\definition{
We call a cocompact free proper group action $G\acts M$ by a (discrete) group $G$ compatible with this structure, and write $G\acts (M\to B, \nabla, g)$, if the Riemannian metric $g$ is $G$-invariant and
there is a group action $G'\acts B$  together with a  surjective group homomorphism $\varphi\colon G\twoheadrightarrow G'$ such that the projection $M\xrightarrow{\pi} B$ is $\varphi$-equivariant\footnote{For all $\gamma\in G$ and $x\in M$ we have $\pi(\gamma x) = \varphi(\gamma)\pi(x)$. In particular, $\ker(\varphi)$ acts on each fibre $F'_b$.}.
}

\example{
The typical example for such a Riemannian fibre bundle with connection and compatible group action is obtained by starting with a compact fibre bundle $F\to M\to B$ where $F$, $M$ and $B$ are connected.
The universal covering $\widetilde{M}$ of $M$ can be considered as a fibre bundle $\widetilde{M}\to\widetilde{B}$ over the universal covering of the base $B$ with some fibres $F'_\bullet$ (in general, these are not the universal coverings of the fibres $F_\bullet$). On the universal coverings, we have the action of $\pi_1(M)$ on $\widetilde{M}$ and the action of $\pi_1(B)$ on $\widetilde{B}$, compare the following diagram:
$$\begin{tikzcd}[ampersand replacement=\&, row sep=10, column sep=20]
\& \pi_1(M) \ar[r, "\varphi"] \ar[d, phantom, "\acts" description] \& \pi_1(B) \ar[d, phantom, "\acts" description]\\
F' \ar[r] \& \widetilde{M} \ar[r]\ar[ddd]\& \widetilde{B} \ar[ddd] \\ \\ \\
F \ar[r] \& M \ar[r] \& B
\end{tikzcd}$$
Here, the long exact sequence of homotopy groups for the fibre bundle $F\to M\to B$, 
$$ \cdots \to \pi_2(B) \to \pi_1(F) \to \pi_1(M) \overset{\varphi}{\twoheadrightarrow} \pi_1(B) \to 0 , $$ 
yields a group homomorphism $\varphi\colon \pi_1(M)\to \pi_1(B)$ that is surjective since $\pi_0(F)$ is trivial. 
The elements in the kernel of $\varphi$ are in the image of $\pi_1(F)\to \pi_1(M)$ and act fibrewise on each fibre $F'_b$ for $b\in \widetilde{B}$ and the projection $\widetilde{M}\to \widetilde{B}$ is $\varphi$-equivariant.
 }

\definition{
Let $(M\to B, \nabla, g)$ be a Riemannian fibre bundle with connection.
For smooth positive functions $s_H,s_V \in \cal C^\infty(M,\R_+)$ we define the Riemannian metric $g^{s_H,s_V}$ on $M$ by 
$$ x\mapsto g^{s_H,s_V}_x =  s_H(x)^2 g_{x,H} + s_V(x)^2 g_{x,V}. $$
In particular, if $s_H\equiv \overline\mu > 0$ and $s_V\equiv \overline\nu > 0$ are constant functions, this defines 
$$g^{\overline\mu,\overline\nu}=g^{s_H,s_V} = \overline \mu ^2 g_V + \overline\nu^2 g_H. $$
}
We use this structure to define a refined version of the spectral density function depending on two parameters in place of the classical parameter $\lambda$.

\definition{
Let $G\acts (M\to B,\nabla,g)$ be a Riemannian fibre bundle with connection and compatible $G$-action. 
Then, using the previous definition, we define the two-parameter spectral density function $\cal G_k(M\to B, \nabla, g)\colon \R_+\times\R_+ \to [0,\infty]$ by
\begin{align*}
    \cal G_k(M\to B, \nabla, g)(\overline\mu,\overline\nu) &=  \tr_{\cal NG} \chi_{[0,1]}(\Delta^k_\up(M, g^{\overline\mu,\overline\nu})) = \cal F_k(M,g^{\overline\mu,\overline\nu})(1).
\end{align*}
We call two such functions $\cal G, \cal G'\colon \R_+\times \R_+\to [0,\infty]$ dilatationally equivalent if there exists a constant $C>0$ such that for all $\overline\mu,\overline\nu\in \R_+$,
$$  \cal G(C^{-1}\overline\mu,C^{-1}\overline\nu) \leq \cal G'(\overline\mu,\overline\nu) \leq \cal G(C\overline\mu, C\overline\nu). $$ 
In this case we write $\cal G\sim \cal G'$.
}

The fact that we chose the value one for the upper end of the interval is not of importance here, in the sense that the dilatational equivalence class of $\cal G$ does not depend on the upper end.

\lemma{
The dilatational equivalence class is independent of the right end chosen for the interval, that is for all $\lambda_0>0$, 
$$ \cal G_k(M\to B, \nabla, g)(\overline\mu,\overline\nu) \sim \left( (\overline\mu,\overline\nu) \mapsto \tr_{\cal NG} \chi_{[0,\lambda_0]}(\Delta^k_\up(M, g^{\overline\mu,\overline\nu}))\right). $$
}
\begin{proof}
This follows directly with constant $C=\sqrt\lambda_0$ since
\begin{align*}
\tr_{\cal NG} \chi_{[0,\lambda_0]}(\Delta^k_\up(M, g^{\overline\mu,\overline\nu}))
&= 
\tr_{\cal NG} \chi_{[0,1]}(\lambda_0^{-1} \Delta^k_\up(M, g^{\overline\mu,\overline\nu}))
\\
&= 
\tr_{\cal NG} \chi_{[0,1]}(\Delta^k_\up(M, \lambda_0 g^{\overline\mu,\overline\nu}))
\\
&= 
\tr_{\cal NG} \chi_{[0,1]}(\Delta^k_\up(M,  g^{\sqrt{\lambda_0}\cdot \overline\mu, \sqrt{\lambda_0} \cdot \overline\nu}))
\\
&= 
\cal G_k(M\to B, g, \nabla)(\sqrt{\lambda_0}\cdot\overline\mu,\sqrt{\lambda_0}\cdot\overline\nu). \qedhere
\end{align*}
\end{proof}

Instead of having two truely independent parameters $\overline{\mu}$ and $\overline{\nu}$, we would like to consider the two parameters as different speeds of scaling the manifold. 
Therefore, we replace these two parameters with two functions, depending on the same variable $\lambda$, governing how fast the fibre respectively the base get scaled as $\lambda\searrow 0$.

\definition{
\label{D_2ParamNSI}
Let $G\acts (M\to B,\nabla,g)$ be a Riemannian fibre bundle with connection and compatible $G$-action. 
Let $\mu,\nu\colon \R_{\geq 0}\to \R_{\geq 0}$ be monotonously increasing continuous functions with $\mu(0)=0=\nu(0)$.
Denoting $\cal G_k = \cal G_k(M\to B, \nabla ,g)$ we define the two-parameter Novikov-Shubin numbers by
\begin{align*}
 \alpha_k(M\to B,\nabla,g)(\mu,\nu) &= \alpha\big(\ \lambda\mapsto  \cal G_k(\mu(\lambda),\nu(\lambda))\ \big) \\
&=\liminf_{\lambda\searrow 0}\frac{\log(\cal G_k(\mu(\lambda),\nu(\lambda)) - b^{(2)}(d_{k+1}))}{\log(\lambda)}. 
\end{align*}
}
Recall here that $b^{(2)}(d_{k+1})$ measures the size of the kernel of $d_{k+1}$ and is metric invariant. Thus, we extend the definition of $\cal G_k$ formally by $\cal G_k(\mu(0),\nu(0))= \cal G_k(0,0)=b^{(2)}(d_{k+1})$. 

\remark{
This two-parameter function generalises the usual spectral density function. Indeed, if\footnote{By abuse of notation, $\lambda$ denotes the function $\id\colon \lambda\mapsto \lambda$ or, more generally, $\lambda^c$ the function $\lambda\mapsto \lambda^c$.} $\mu=\nu=\lambda$, then $g^{\lambda,\lambda} = \lambda^2 g =  g_\lambda $ independently of the connection $\nabla$ chosen. Hence,
$$ \cal G_k(M\to B, \nabla, g)(\lambda,\lambda) = \tr_{\cal NG} \chi_{[0,1]}(\Delta^k_\up(M, g_\lambda)) =  \cal F_k(M,g)(\lambda)  $$
is the classical spectral density function of $(M,g)$ and therefore 
$$ \alpha_k(M\to B, g, \nabla)(\lambda,\lambda) = \alpha_k(M) $$
recovers the Novikov-Shubin invariants.
}

\example{
In the simplest case of a product manifold $(M,g) = (F,g_F)\times (B,g_B)$ with the canonical connection $TM \cong_\nabla TF\perp TB$, for $\mu,\nu>0$ we have
\begin{align*}
    \cal G_k(F\times B, \nabla, g)(\mu,\nu) 
    &= \cal F_k((F,\nu^2 g_F)\times (B,\mu^2 g_B))(1).
\end{align*}
By \cite[Cor.~2.44]{LckL2}, it is therefore dilatationally equivalent to 
\begin{align*}
   \cal G_k(F\times B, \nabla, g)(\mu,\nu)  &\sim \sum_{p+q=k} \cal F_p((F, \nu^2 g_F))(1)\cdot \cal F_q((B,\mu^2 g_B))(1)
    \\
    &= \sum_{p+q=k} \cal F_p(F)(\nu) \cdot \cal F_q(B)(\mu).
\end{align*}
If $\mu = \lambda^r$ and $\nu = \lambda^s$, we can consider a limit as $\lambda\searrow 0$ in the spirit of the Novikov-Shubin invariants. 
We assume that all $L^2$-Betti numbers in this example vanish\footnote{This is not necessary but reduces the length of notation for this example considerably. One can proceed just as in cited source by W.~Lück even if the $L^2$-Betti numbers do not vanish.}. Following the computation in W.~Lück's book~\cite[Thm.~2.55~(3)]{LckL2},
\begin{align*}
    \alpha_k(F \times B, \nabla, g)(\mu,\nu) 
    &= \liminf_{\lambda\searrow 0} \frac{\log(\cal G_k(F\times B, \nabla,g)(\lambda^r,\lambda^s))}{\log(\lambda)}
    \\
    &= \liminf_{\lambda\searrow 0} \frac{\log(\cal F_k((F,\lambda^{2s} g_F)\times (B,\lambda^{2r} g_B), \nabla,g)(1))}{\log(\lambda)}
    \\
    &= \min_{0\leq p\leq k}\left\{ \begin{matrix} \alpha\left(\cal F_p(F)(\lambda^s) \cdot \cal F_{k-p}(B)(\lambda^r)\right),\\ \alpha\left(\cal F_{p+1}(F)(\lambda^s) \cdot \cal F_{k-p}(B)(\lambda^r)\right) \end{matrix} \right\}
    \\
    &= \min_{0\leq p\leq k}\left\{ \begin{matrix} \alpha\left(\cal F_p(F)(\lambda^s)\right) + \alpha\left(\cal F_{k-p}(B)(\lambda^r)\right),\\ \alpha\left(\cal F_{p+1}(F)(\lambda^s)\right) + \alpha\left(\cal F_{k-p}(B)(\lambda^r)\right) \end{matrix} \right\}
    \\
    &= \min_{0\leq p\leq k}\left\{ \begin{matrix} s\cdot\alpha_p(F) + r\cdot\alpha_{k-p}(B),\\ s\cdot\alpha_{p+1}(F) + r\cdot\alpha_{k-p}(B) \end{matrix} \right\}.
\end{align*}
In this case, we see the contributions from the base and fibre are scaled according to the chosen functions $\mu(\lambda)=\lambda^r$ and $\nu(\lambda)=\lambda^s$ as $\lambda\searrow 0$.
}

%% file: NearCohom.tex
When studying the classical spectral density functions and Novikov-Shubin invariants, a useful approach in many cases is by using the notion of near cohomology, as given by M.~Gromov and M.~A.~Shubin in their paper~\cite{GroShu}.
We show that a similar approach can still be used in this two-parameter case.

Decomposing the tangent bundle $TM\xrightarrow{\pi} M$  as $TM \cong VM\oplus HM$ into a vertical and a horizontal subbundle gives us a diagram
$$
\begin{tikzcd}[ampersand replacement = \&]
TF_\bullet \ar[d] \ar[r] 
\& 
T_\bullet F_\bullet \ar[dr]
\& 
\begin{matrix}  VM \oplus HM \\ \cong \\ TM \end{matrix} \ar[d] 
\ar[l, "\cong", no head, rounded corners, to path ={[pos=.15]([yshift=.45cm]\tikztostart.west) -| ([yshift=-3pt]\tikztotarget.north)\tikztonodes -- ([yshift=-3pt]\tikztotarget.north)}] 
\ar[r, "\cong"', no head, rounded corners, to path ={[pos=.15]([yshift=.45cm]\tikztostart.east) -| ([yshift=-3pt]\tikztotarget.north)\tikztonodes -- ([yshift=-3pt]\tikztotarget.north)}]
\arrow[rr, "d\pi"',rounded corners,
  to path={[pos=0.25]  
  (\tikztostart.north)
    -| ([yshift=.5cm]\tikztostart.north)
    -| (\tikztotarget.north)\tikztonodes
    -- (\tikztotarget.north)} ]
\&
\pi^*TB \ar[dl] \ar[r, "d\pi"]
\&  
TB \ar[d]
\\
F_\bullet \ar[rr]
\& \&
M \ar[rr, "\pi"]
\& \&
B.
\end{tikzcd}
$$
Given any vector field $X\in \Gamma(TM)$, we can decompose 
$$ X = Y + Z, \qquad \text{with} \qquad Y\in \Gamma(\pi^*TB), \quad Z\in \Gamma(T_\bullet F_\bullet)$$
into a horizontal component $Y$ and a vertical component $Z$.

We call a vector field $Y\in \Gamma(\pi^*TB)$ basic, 
if there exists a vector field $\overline{Y}\in \Gamma(TB)$ such that $Y$ is $\pi$-related to $\overline{Y}$, that is, the following diagram commutes (compare, for example, Besse~\cite[Ch. 9]{Bes}):
$$\begin{tikzcd}[ampersand replacement = \&]
 TM \ar[r, "d\pi"] 
 \&
 TB \\
 M \ar[u, "Y"] \ar[r, "\pi"']
 \&
 B \ar[u, "\overline{Y}"']
\end{tikzcd}$$  
We call $Y$ the lift of $\overline{Y}$. 
For every $U\in \Gamma(TB)$ there exists a unique such lift $\widetilde U\in \Gamma(\pi^* TB)$.
We denote by $\Gamma_b(HM) \subset \Gamma(\pi^*TB)$ the set of basic vector fields. 
Then $\Gamma_b(HM)$ spans $\Gamma(\pi^*TB)$ as a $\cal C^\infty(M)$-module, so every horizontal vector field $Y\in \Gamma(HM)\cong \Gamma(\pi^*(TB))$ can be written as
$$ Y = \sum_{i\in I} f_i\cdot \widetilde{U_i}$$
for smooth functions $f_i \in \cal C^\infty(M)$ and $U_i\in \Gamma(TB)$.

\lemma{
	Let $Z,Z'\in \Gamma(VM)$ be vertical vector fields and $Y\in \Gamma_b(HM)$ a basic horizontal vector field.
	Then 
	\begin{enumerate}
		\item $[Z, Z'] \in \Gamma(VM)$,
		\item $[Y, Z] \in \Gamma(VM)$.
	\end{enumerate}
}
\begin{proof}
Recall that $VM = \ker(d\pi)$, hence $Z\sim_\pi 0$ and $Z\sim_\pi 0$ where $0\in \Gamma(TB)$ denotes the zero section.
By definition, $Y\sim_\pi \overline{Y}$ for some $\overline{Y}\in TB$. Therefore,
\begin{align*}
 d\pi [Z, Z'] 
 &= \widetilde{[0,0]_B} = 0, \qquad
 d\pi [Y,Z] 
 = \widetilde{[\overline{Y}, 0]_B} = 0,
\end{align*}
and the claim follows.
\end{proof}

Looking at the de Rham complex $\Omega^\bullet(M)$, it can be decomposed using the fibre bundle structure.

\theorem{
\label{T_splittingOmegaFMB}
Let $F_\bullet\to M\to B$ be a fibre bundle, then there is an isomorphism
$$ \Omega^k(M) \xrightarrow[\ \cong\ ]{\ \Phi\ } \bigoplus_{p+q=k} \Omega^p(B, \{\Omega^q(F_b)\}_{b\in B}), $$
identifying forms on $M$ and forms on $B$ with values in the system of forms on the fibres $\{F_b\}_{b\in B}$.\footnote{The right-hand-side is understood in the sense of A.~Fomenko and D.~Fuchs~\cite[Lec.~22.2]{FomenkoFuchs}.}
}

\begin{proof}
Using that $TM\cong VM\oplus HM$, we decompose $X\in \Gamma(TM)$ as $X = Y+Z$ with $Y\in \Gamma(HM)$ and $Z\in \Gamma(VM)$. 
Given $U_1,\dots, U_p\in \Gamma(TB)$ with basic lifts $\widetilde{U_1},\dots, \widetilde{U_p}\in \Gamma_b(HM)$ and $Z_{p+1},\dots, Z_k\in T_\bullet F_\bullet \cong VM$,
for a $k$-form $\omega\in \Omega^k(M)$ we define 
$$ \Phi(\omega) = \sum_{p+q=k} (\Phi(\omega))_{p,q} $$
where the $(p,q)$-summand $(\Phi(\omega))_{p,q} \in \Omega^p(B,\{\Omega^q(F_\bullet)\})$ is given by
$$(\Phi(\omega))_{p,q}(U_1,\dots,U_p)(Z_{p+1},\dots,Z_k) = \omega\left(\widetilde{U_1},\dots,\widetilde{U_p},Z_{p+1},\dots, Z_k \right). $$

Decomposing $X_\bullet\in \Gamma(TM)$ as $X_\bullet=Y_\bullet+Z_\bullet$ with $Y_\bullet\in \Gamma(HM)$ and $Z_\bullet\in \Gamma(VM)$ as before, we construct the inverse 
$$\Psi\colon \bigoplus_{p+q=n} \Omega^p(B,\{\Omega^q(F_\bullet)\})\to\Omega^n(M)$$
 to this map, starting with $\alpha\in\Omega^p(B,\{\Omega^q(F_\bullet)\})$ by 
\begin{align*}
 \Psi(\alpha)(X_1,\dots, X_k) 
 &= \frac{1}{p!q!}\sum_{\sigma\in \cal S_k}\mathrm{sgn}(\sigma) \left(\pi^*\alpha\right)( Y_{\sigma(1)},\dots, Y_{\sigma(p)})( Z_{\sigma(p+1)},\dots, Z_{\sigma(k)}),
\end{align*}
where $\cal S_k$ is the set of permutations of the first $k$ integers, $\{1,\dots,k\}$, and $\mathrm{sgn}$ the sign of the permutation. 
This is then extended linearly to the direct sum.

We check that $\Psi$ and $\Phi$ are indeed inverses to each other. 
With the notation above, for a summand $\alpha\in\Omega^p(B,\{\Omega^q(F_\bullet)\})$,
\begin{align*}
\Phi\Psi(\alpha)&(U_1,\dots, U_p)(Z_{p+1},\dots, Z_k) 
\\
&= \Psi(\alpha)\left(\widetilde{U_1},\dots, \widetilde{U_p},Z_{p+1},\dots, Z_k\right)
\\
&= 
\frac{1}{p!q!}\sum_{\sigma\in S_k}\mathrm{sgn}(\sigma) \left(\pi^*\alpha\right)\left(\widetilde{U_{\sigma(1)}},\dots, \widetilde{U_{\sigma(p)}}\right)( Z_{\sigma(p+1)},\dots, Z_{\sigma(k)})
\\
&=
\left(\pi^*\alpha\right)\left(\widetilde{U_1}, \dots, \widetilde{U_p}\right)(Z_{p+1},\dots,Z_n)
\\
&=
\alpha(U_1,\dots, U_p)(Z_{p+1},\dots,Z_k),
\end{align*}
where in the third equality $U_l=0$ for $l>p$ and $Z_l=0$ for $l\leq p$, so that after reordering the arguments, each summand appears $p!q!$ times with $+$ sign.

In the other direction, writing $X_\bullet = Y_\bullet+Z_\bullet\in \Gamma(HM)\oplus \Gamma(VM) \cong \Gamma(TM)$ as before,
\begin{align*}
\Psi\Phi(\omega)&(X_1,\dots, X_k) \\
&= \sum_{p+q=k}\frac{1}{p!q!}\sum_{\sigma\in S_k}\mathrm{sgn}(\sigma) (\pi^*\Phi(\omega))(Y_{\sigma(1)},\dots, Y_{\sigma(p)})(Z_{\sigma(p+1)},\dots, Z_{\sigma(k)}),
\end{align*}
where pointwise for $x\in M$ with $b=\pi(x)$, 
\begin{align*}
(\pi^*\Phi(\omega))_x&(Y_{\sigma(1)}(x),\dots, Y_{\sigma(p)}(x)) (Z_{\sigma(p+1)}(x),\dots, Z_{\sigma(k)}(x))
\\
&= \Phi(\omega)_b(Y_{\sigma(1)}(x),\dots, Y_{\sigma(p)}(x)) (Z_{\sigma(p+1)}(x),\dots, Z_{\sigma(k)}(x)) 
\\
&= \omega_x(A_{\sigma(1)}(x),\dots, A_{\sigma(p)}(x), Z_{\sigma(p+1)}(x),\dots, Z_{\sigma(k)}(x))
\\
&= \omega_x({Y_{\sigma(1)}(x)},\dots, {Y_{\sigma(p)}(x)}, Z_{\sigma(p+1)}(x),\dots, Z_{\sigma(k)}(x))
\end{align*}
where $A_i$ is some basic horizontal vector field with $A_i(x) = Y_i(x)$. 
Therefore, 
\begin{align*}
\Psi\Phi(\omega)(X_1,\dots, X_k) 
&=
\sum_{p+q=k}\frac{1}{p!q!}\sum_{\sigma\in S_n}\mathrm{sgn}(\sigma) \omega(Y_{\sigma(1)},\dots, Y_{\sigma(p)},Z_{\sigma(p+1)},\dots, Z_{\sigma(k)})
\\
&= \sum_{(\Xi_1,\dots,\Xi_k) \in \{Y_1,Z_1\}\times\cdots\times \{Y_k,Z_k\}} \omega(\Xi_1, \dots, \Xi_k) 
\\
&= \omega(X_1,\dots, X_k),
\end{align*}
where we use in the second equality that $\omega$ is antisymmetric and that after ordering each summand appears $p!q!$ times, where $p$ is the number of $Y_\bullet$s and $q$ the number of $Z_\bullet$s chosen. The last equality then follows by linearity of $\omega$.
\end{proof}

We can now look at the differential $d\colon \Omega^k(M)\to \Omega^{k+1}(M)$ under this decomposition.

\lemma{
	Under the decomposition $\Phi$ of $\Omega^\bullet(M)$, the de Rham differential splits into three summands, 
	$ d \cong d^{0,1} + d^{1,0}  + d^{2,-1}, $
	where 
	$$d^{i,1-i}\colon \Omega^p(B,\{\Omega^q(F_\bullet)\}) \to \Omega^{p+i}(B,\{\Omega^{q+1-i}(F_\bullet)\}).$$
}
\begin{proof}
By Cartan's formula, for $\omega\in \Omega^k(M)$ and $X_0,\dots, X_{k}\in \Gamma(TM)$, the de Rham differential of $\omega$ evaluated on the $X_\bullet$s is given by
\begin{align*}
d(\omega)(X_1,\dots, X_{k+1})
&=
\sum_{i=0}^k (-1)^i X_i ( \omega( X_0,\overset{\widehat X_i}{\dots}, X_k) ) \\
& + \sum_{0\leq i< j\leq k} (-1)^{i+j} \omega([X_i, X_j],X_0,\overset{\widehat X_i, \widehat X_j}{\dots}, X_k).
\end{align*}
We denote by $[X_i, X_j]_H$ respectively $[X_i, X_j]_V$ the projection of $[X_i, X_j]$ to $\Gamma(HM)$ respectively $\Gamma(VM)$.
Given $\alpha\in \Omega^p(B, \{\Omega^q F_\bullet\})$, we compute 
$\Phi d\Psi \alpha$ by looking at the $(r,s)$-component $(\Phi d \Psi \alpha)_{r,s}\in \Omega^r(B, \{\Omega^s F_\bullet\})$ (with $r+s=k+1)$.

For this, let $U_1,\dots, U_r\in \Gamma(TB)$ and $Z_{r+1},\dots, Z_{k+1}\in \Gamma(T_\bullet F_\bullet)$, then
\begin{align*}
(\Phi d \Psi \alpha)_{r,s}&(U_1,\dots, U_r)(Z_{r+1},\dots, Z_{k+1})
\\
&=(d\Psi\alpha)\left(\widetilde{U_1},\dots,\widetilde{U_r},Z_{r+1},\dots, Z_{k+1}\right)
\\
&= \sum_{1\leq i\leq r} (-1)^{i+1} \widetilde{U_i}\left(\Psi\alpha(\widetilde{U_1},\overset{\widehat{{U_i}}}{\dots},\widetilde{U_r},Z_{r+1},\dots,Z_{k+1})\right) 
\\
&\quad + \sum_{r+1\leq i\leq k+1} (-1)^{i+1} Z_i\left(\Psi\alpha(\widetilde{Y_1},\dots,\widetilde{U_r},Z_{r+1},\overset{\widehat{Z_i}}{\dots},Z_{k+1})\right) 
\\
&\quad + \sum_{1\leq i<j\leq r} (-1)^{i+j+1} \Psi\alpha\left([\widetilde{U_i},\widetilde{U_j}],\widetilde{U_1},\overset{\widehat{{U_i}},\widehat{U_j}}{\dots},\widetilde{U_r}, Z_{r+1},\dots,Z_{k+1}\right) 
\\
&\quad + \sum_{1\leq i \leq r < j\leq k+1} (-1)^{i+j+1} \Psi\alpha\left([\widetilde{U_i},Z_j],\widetilde{U_1},\overset{\widehat{U_i}}{\dots},\widetilde{U_r}, Z_{r+1},\overset{\widehat{Z_j}}{\dots},Z_{k+1}\right) 
\\
&\quad + \sum_{r+1\leq i<j\leq k+1} (-1)^{i+j+1} \Psi\alpha\left([Z_i,Z_j],\widetilde{U_1},\dots,\widetilde{U_r}, Z_{r+1},\overset{\widehat{Z_i},\widehat{Z_j}}{\dots},Z_{k+1}\right).
\end{align*} 

By definition, $\Psi(\alpha)\neq 0$ only if $p$ of the arguments have non-zero components in $\Gamma(HM)$ and $q$ of the arguments have non-zero components in $\Gamma(VM)$.
Recall that $[Z,Z'],[\widetilde{U},Z]\in \Gamma(VM)$ for all $Z,Z'\in \Gamma(VM)$ and $U\in \Gamma(TB)$.
Therefore, the operator $\Phi d\Psi \alpha$ decomposes into the following three summands.
\begin{enumerate}
\item
The first summand keeps the base-degree fixed and increases the fibre-degree by one. 
It is given for $\alpha\in \Omega^p(B, \{\Omega^q F_\bullet\})$ by
\begin{align*}
(\Phi d \Psi \alpha)_{p,q+1}&(U_1,\dots, U_p)(Z_{p+1}, \dots, Z_{k+1})\\
&= \sum_{p+1\leq i\leq k+1} (-1)^{i+1} Z_i(\Psi\alpha(\widetilde{U_1},\dots,\widetilde{U_p},Z_{p+1},\overset{\widehat{Z_i}}{\dots},Z_{k+1})) \\
&\quad + \sum_{p+1\leq i<j\leq k+1} (-1)^{i+j+1-p} \Psi\alpha(\widetilde{U_1},\dots,\widetilde{U_p}, [Z_i,Z_j], Z_{p+1},\overset{\widehat{Z_i},\widehat{Z_j}}{\dots},Z_{k+1})
\\
&=
\sum_{p+1\leq i\leq k+1} (-1)^{i+1} Z_i(\alpha(U_1,\dots, U_p))(Z_{p+1},\overset{\widehat{Z_i}}{\dots}, Z_{k+1})
\\
&\quad + \sum_{p+1\leq i<j \leq k+1} (-1)^{i+j+1-p} \alpha(U_1,\dots, U_p)([Z_i, Z_j], Z_{p+1}, \overset{\widehat{Z_i},\widehat{Z_j}}{\dots}, Z_{k+1}).
\end{align*}
\item
The second summand increases the base-degree by one and keeps the fibre-degree fixed. It is given by
\begin{align*}
(\Phi d \Psi \alpha)_{p+1,q}&(U_1,\dots, U_{p+1})(Z_{p+2}, \dots, Z_{k+1})\\
&= \sum_{1\leq i\leq p+1} (-1)^{i+1} \widetilde{U_i}(\Psi\alpha(\widetilde{U_1},\overset{\widehat{U_i}}{\dots},\widetilde{U_{p+1}},Z_{p+2},\dots,Z_{k+1})) \\
&\quad + \sum_{1\leq i<j\leq p+1} (-1)^{i+j+1} \Psi\alpha([\widetilde{U_i},\widetilde{U_j}]_H,\widetilde{U_1},\overset{\widehat{U_i},\widehat{U_j}}{\dots},\widetilde{U_{p+1}}, Z_{p+2},\dots,Z_{k+1}) \\
&\quad + \sum_{1\leq i \leq p+1 < j\leq k+1} (-1)^{i+j+1-p} \Psi\alpha(\widetilde{U_1},\overset{\widehat{U_i}}{\dots}, \widetilde{U_{p+1}}, [\widetilde{U_i},Z_j], Z_{p+2},\overset{\widehat{Z_j}}{\dots},Z_{k+1})
\\
&= 
\sum_{1\leq i\leq p+1} (-1)^{i+1} {\widetilde U_i}(\alpha({U_1},\overset{\widehat{U_i}}{\dots},{U_{p+1}})(Z_{p+2},\dots,Z_{k+1})) \\
&\quad + \sum_{1\leq i<j\leq p+1} (-1)^{i+j+1} \alpha([{U_i},{U_j}],{U_1},\overset{\widehat{U_i},\widehat{U_j}}{\dots},{U_{p+1}})( Z_{p+2},\dots,Z_{k+1}) \\
&\quad + \sum_{1\leq i \leq p+1 < j\leq k+1} (-1)^{i+j+1-p} \alpha({U_1},\overset{\widehat{U_i}}{\dots}, {U_{p+1}})( [\widetilde{U_i},Z_j], Z_{p+2},\overset{\widehat{Z_j}}{\dots},Z_{k+1}).
\end{align*}
\item 
The third summand increases the base-degree by two and decreases the fibre-degree by one. It is given by
\begin{align*}
(\Phi d \Psi \alpha)_{p+2,q-1}&(U_1,\dots, U_{p+2})(Z_{p+3}, \dots, Z_{k+1})\\
&= \sum_{1\leq i<j\leq p+2} (-1)^{i+j+1-p} \Psi\alpha(\widetilde{U_1},\overset{\widehat{U_i},\widehat{U_j}}{\dots}, \widetilde{U_{p+2}}, [\widetilde{U_i},\widetilde{U_j}]_V, Z_{p+3},\dots,Z_{k+1})
\\
&= \sum_{1\leq i<j\leq p+2} (-1)^{i+j+1-p} \alpha({U_1},\overset{\widehat{U_i},\widehat{U_j}}{\dots}, {U_{p+2}})([\widetilde{U_i},\widetilde{U_j}]_V, Z_{p+3},\dots,Z_{k+1}).
\end{align*}
\end{enumerate}

The claim follows with the maps defined for $\alpha\in \Omega^p(B,\{\Omega^q(F_\bullet)\})$ by
\begin{align*}
d^{0,1}(\alpha) &= (\Phi d\Psi\alpha)_{p,q+1},
\\
d^{1,0}(\alpha) &= (\Phi d\Psi\alpha)_{p+1,q},
\\
d^{2,-1}(\alpha) &= (\Phi d\Psi\alpha)_{p+2,q-1}.
\end{align*}
and $d=d^{0,1} + d^{1,0} + d^{2,-1}$ extended linearly to $\bigoplus_{p+q=k}\Omega^p(B,\{\Omega^q(F_\bullet)\})$.
\end{proof}

Denote $E_0^{p,q}  = \Omega^{p}(B, \{\Omega^q(F_\bullet)\})$, then we can visualise this decomposition as a $\Z^2$-graded complex.\footnote{This is not a double complex in general as there is the diagonal $d^{2,-1}$-map. If we can choose a flat connection on $M$, then $d^{2,-1}$ vanishes and this is a true double complex. In terms of objects, this may be viewed as the zeroth page of the Serre spectral sequence of $F_\bullet\to M\to B$.}
An excerpt of this is pictured below, with the maps $d^{i,1-i}$ only drawn at $E_0^{p,q}$ and as dashed arrows at their images.
As usual, the parts appearing in $\Omega^k(M)$ align along the antidiagonal $p+q=k$ in the diagram.
$$
\fbox{
\begin{tikzcd}[ampersand replacement = \&]
E_0^{p,q+2}
\&
E_0^{p+1,q+2}
\&
E_0^{p+2,q+2}
\&
E_0^{p+3,q+2}
\&
E_0^{p+4,q+2}
\\
E_0^{p,q+1}
\ar[u, dashed] \ar[r, dashed] \ar[rrd, dashed]
\&
E_0^{p+1,q+1}
\&
E_0^{p+2,q+1}
\&
E_0^{p+3,q+1}
\&
E_0^{p+4,q+1}
\\
E_0^{p,q}
\ar[u, "d^{0,1}" description] \ar[r, "d^{1,0}" description] \ar[rrd, "d^{2,-1}" description]
\&
E_0^{p+1,q}
\ar[u, dashed] \ar[r, dashed] \ar[rrd, dashed]
\&
E_0^{p+2,q}
\&
E_0^{p+3,q}
\&
E_0^{p+4,q}
\\
E_0^{p,q-1}
\&
E_0^{p+1,q-1}
\&
E_0^{p+2,q-1}
\ar[u, dashed] \ar[r, dashed] \ar[rrd, dashed]
\&
E_0^{p+3,q-1}
\&
E_0^{p+4,q-1}
\\
E_0^{p,q-2}
\&
E_0^{p+1,q-2}
\&
E_0^{p+2,q-2}
\&
E_0^{p+3,q-2}
\&
E_0^{p+4,q-2}
\end{tikzcd}
}$$

Since $d = d^{0,1}+d^{1,0}+d^{2,-1}$ is a differential, that is, $d^2 = 0$, we obtain immediately that
\begin{alignat*}{2}
 &0 = (d^{0,1})^2, 
 &&0 = d^{0,1} d^{1,0} + d^{1,0} d^{0,1},  \\
&0 = d^{0,1}d^{2,-1}  + (d^{1,0})^2 + d^{2,-1}d^{0,1}, \qquad  \qquad
 &&0 = d^{1,0}d^{2,-1} + d^{2,-1}d^{1,0}, \\
 &0 =  (d^{2,-1})^2. &&
\end{alignat*}

Note that $d^{1,0}$ is not a differential in general. 
Leaving out the terms that cancel due to the usual alternating sign\footnote{Coming from leaving out two arguments in two different orders.}, a direct computation shows that for $\alpha\in E_0^{p,q}$, $U_1,\dots,U_{p+2}\in \Gamma(TB)$ and $Z_{p+3},\dots, Z_{k+2}\in \Gamma(T_\bullet F_\bullet)$:
\begin{align*}
(d^{1,0})^2(\alpha)&(U_1,\dots, U_{p+3})(Z_{p+3},\dots, Z_{k+2})
\\
&=
\sum_{1\leq i < j\leq p+2} (-1)^{i+j} (\widetilde U_j \widetilde U_i - \widetilde U_i \widetilde U_j)(\alpha(U_1,\overset{\widehat U_i, \widehat U_j}{\dots}, U_{p+2})(Z_{p+3},\dots, Z_{k+2}))
\\
&\quad +
\sum_{1\leq i < j\leq p+2} (-1)^{i+j} \widetilde{[U_i,U_j]_B} (\alpha(U_1,\overset{\widehat U_i, \widehat U_j}{\dots}, U_{p+2})(Z_{p+3},\dots, Z_{k+2}))
\\
&\quad +
\sum_{1\leq i< j\leq p+2 <l\leq k+2} (-1)^{i+j+l+p}\alpha(U_1,\overset{\widehat U_i, \widehat U_j}{\dots}, U_{p+2})([\widetilde U_i, [\widetilde U_j, Z_l]], Z_{p+3},\dots, Z_{k+2}).
\end{align*}
For the first two terms we have
$$ [\widetilde U_i, \widetilde U_j] - \widetilde{[U_i, U_j]_B}= [\widetilde U_i, \widetilde U_j]_V $$
 and since by the Jacobi identity $[\widetilde U_i, [\widetilde U_j, Z_l]] - [\widetilde U_j, [\widetilde U_i, Z_l]] = -[[\widetilde U_i, \widetilde U_j], Z_l]$
this precisely cancels out the terms that survive in $d^{0,1}d^{2,-1} + d^{2,-1}d^{0,1}$,
\begin{align*}
(d^{0,1}&d^{2,-1} + d^{2,-1}d^{0,1})(\alpha)(U_1,\dots, U_{p+3})(Z_{p+3},\dots, Z_{k+2})
\\
&=
\sum_{1\leq i < j\leq p+2} (-1)^{i+j} [\widetilde U_i, \widetilde U_j]_V (\alpha(U_1,\overset{\widehat U_i, \widehat U_j}{\dots}, U_{p+2})(Z_{p+3},\dots, Z_{k+2}))
\\
&\quad +
\sum_{1\leq i< j\leq p+2 <l\leq k+2} (-1)^{i+j+l+p+1}\alpha(U_1,\overset{\widehat U_i, \widehat U_j}{\dots}, U_{p+2})([[\widetilde U_i, \widetilde U_j], Z_l], Z_{p+3},\dots, Z_{k+2}).
\end{align*}

The inner product on $\Omega^k(M)$ (coming from the Riemannian metric) induces inner products on the decomposition.
For $\alpha,\alpha'\in \Omega^p(B, \{\Omega^q(F_\bullet)\})$, 
$$ \langle \alpha,\alpha' \rangle_{g, (p,q)} = \langle \Psi\alpha, \Psi\alpha'\rangle_{(M,g)} = \int_{(M,g)} \Psi\alpha\wedge\ast\Psi\alpha', $$
whereas the different direct summands are mutually orthogonal to each other since the decomposition of $TM$ into $VM$ and $HM$ is orthogonal. 
This implies for $\omega\in \Omega^k(M)$ with $\Phi(\omega)=\alpha = \sum_{p+q=k}\alpha_{p,q} \in \bigoplus_{p+q=k} E_0^{p,q}$, 
$$\|\omega\|_g^2 = \|\alpha\|^2_g = \sum_{p+q=k}\|\alpha_{p,q}\|^2_g = \sum_{p+q=k} \langle \alpha_{p,q}, \alpha_{p,q} \rangle_{g,(p,q)} . $$ 
When changing the metric from $g$ to $g^{\mu,\nu}$ on $M$, the length of a vertical tangent vector $v\in V_xM$ changes by a factor $\nu$ as
$$ \|v\|_{g^{\mu,\nu}}^2 = \nu^2 g_x^V(v,v) = (\nu \|v\|_{g})^2 $$
and on horizontal tangent vectors $h\in H_xM$ by $\|h\|_{g^{\mu,\nu}}^2 = (\mu \|h\|_{g})^2$.
Denote by $\omega_g$ the volume form on $(M,g)$ and by $\omega_{g^{\mu,\nu}}$ the volume form on $(M, g^{\mu,\nu})$. 
By the observation above,
$$ \omega_g = \mu^{-\dim(B)} \nu^{-\dim(F)}\cdot  \omega_{g^{\mu,\nu}}. $$

The Hodge $*$-operators $\ast_g, \ast_{g^{\mu,\nu}}$ map $\Omega^k(M)\to \Omega^{n-k}(M)$ and preserve the decomposition as
$$ \ast\colon \Omega^{p}(B, \{\Omega^{q}(F_\bullet)\}) \to \Omega^{\dim(B)-p}(B, \{\Omega^{\dim(F)-q}(F_\bullet)\}). $$
Since on $\Omega^{p}(B, \{\Omega^{q}(F_\bullet)\})$,
\begin{align*}
\langle \alpha, \beta\rangle_{g^{\mu,\nu}}
&= \int_{(M,g^{\mu,\nu})} \alpha\wedge \ast_{g^{\mu,\nu}} \beta
= \mu^{-\dim(B)} \nu^{-\dim(F)}  \int_{(M,g)} \alpha\wedge \ast_{g^{\mu,\nu}} \beta
\\
&= \mu^{-\dim(B)} \nu^{-\dim(F)} \mu^{\dim(B)-2p} \nu^{\dim(F)-2q}\cdot \int_{(M,g)} \alpha\wedge \ast_{g} \beta
\\
&=\mu^{-2p} \nu^{-2q} \langle \alpha, \beta\rangle_{g},
\end{align*}
the scalar product changes by a factor $\mu^{-2p} \nu^{-2q}$.

This allows us to define the two-parameter Novikov-Shubin numbers via the near cohomology cones of the decomposed complex.
Since the near cohomology cone satisfies
\begin{align*}
&C_{\lambda_0}^k(M,g^{\mu,\nu}) = \left\{ \omega\in \Omega^k(M) \cap \ker(d)^\perp \:\middle|\: \|d\omega\|_{g^{\mu,\nu}} \leq \lambda_0 \|\omega\|_{g^{\mu,\nu}}  \right\}
\\
& \cong \left\{ \! \alpha \! \in \!\! \left(\bigoplus_{p+q=k} E_0^{p,q}\right) \! \cap  \Phi\!\left(\ker(d)^\perp\right) \:\middle|\: \sum_{r+s=k+1} \!\!\!\! \mu^{-r}\nu^{-s}\|(d\alpha)_{r,s}\|_g \leq \lambda_0 \!\!\! \sum_{p+q=k} \!\!\! \mu^{-p}\nu^{-q} \|\alpha_{p,q}\|_g\right\}\!,
\end{align*}
we can define $\cal G_k(M\to B, \nabla, g)$ in terms of this near cohomology cone with $\lambda_0=1$ as follows.
\corollary{
In the notation as above,
$$ \cal G_k(M\to B, \nabla, g)(\mu,\nu) = \sup_{L }\dim_{\cal NG} L, $$
where the supremum runs over all closed linear subspaces $L$ of $C_1^k(M, g^{\mu,\nu})$.
}
\begin{proof}
This follows immediately since 
$ \cal G_k(M\to B, \nabla, g)(\mu,\nu) = \cal F_k(M, g^{\mu,\nu})(1)$.
\end{proof}

%% file: Invariance.tex
In this section we show multiple invariance properties of the two-parameter Novikov-Shubin numbers. 
We show that for a fibre bundle $M\to B$ and a fixed connection $\nabla$, the dilatational equivalence class of the underlying spectral density functions is independent of the $\nabla$-compatible Riemannian metric $g$ on $M$.
Then we show that the spectral density functions are dilatationally equivalent for two bundles $M\to B$ and $M'\to B$ if there exists a certain type of $\nabla$-compatible fibre homotopy equivalence. 
We also show that the dilatational equivalence class of the spectral density functions is independent of the connection $\nabla$ if we restrict them to the parameter subspace $\{\nu\leq \mu\}$, where the fibre is scaled at least as fast as the base.
In particular, the two-parameter Novikov-Shubin numbers are invariant under these operations.

%% file: FixedConnMetricInvar.tex
From the definition in terms of near cohomology cones, we can derive that the dilatational equivalence class of $\cal G_k(M\to B, \nabla, g)$ for a fixed connection $\nabla$ does not depend on the metric $g$.

\theorem{
\label{Thm_2NSI_Metric_fixed_conn_invar}
    Let $G\acts (M\to B,\nabla,g)$ be a fibre bundle with fixed connection $\nabla$ and compatible cocompact free proper group action by a group $G$. Then for $0\leq k\leq \dim(M)$ the dilatational equivalence class of
    $$ \cal G_k(M\to B, \nabla) = \cal G_k(M\to B, \nabla, g)$$
    does not depend on the choice of $G$-invariant $\nabla$-compatible Riemannian metric $g$. 
}
\begin{proof}
On a compact manifold $\overline M$, any two Riemannian metrics $\overline g,\overline{g}'$ are quasi-equivalent, that is there exists $K\geq 1$ such that
$ K^{-1} \overline{g} \leq \overline{g}' \leq K \overline{g}. $
By $G$-invariance of the Riemannian metrics and cocompactness of the action $G\acts M$, this is true for any two choices of $G$-invariant Riemannian metrics $g,g'$ on $M$.
Restricting to the subbundles $V^*M$ and $H^*M$ of $T^*M$, this inequality holds also for the vertical and horizontal parts individually. 
After rescaling, it follows that there is $K>0$ such that for all $\mu,\nu>0$,
$$ K^{-1} g^{\mu,\nu} \leq (g')^{\mu,\nu} \leq K g^{\mu,\nu}. $$
If $\omega\in C_\lambda^k(M,(g')^{\nu,\mu})$, then 
\begin{align*}
     K^{-2(k+1)} \|d\omega\|^2_{g^{\mu,\nu}} &= \|d\omega\|^2_{K^{-1} g^{\mu,\nu}} \leq \|d\omega\|^2_{g'^{\mu,\nu}} 
     \\
     &\leq \lambda^2 \|\omega\|_{g'^{\mu,\nu}}^2 \leq \lambda^2 \|\omega\|^2_{K g^{\mu,\nu}} = K^{2k} \lambda^2 \|\omega\|^2_{g^{\mu,\nu}}.
\end{align*}
This implies that $\omega\in C^k_{K^{2k+1} \lambda}(M, g^{\mu,\nu})= C^k_{\lambda}(M, Kg^{\mu,\nu}).$
We can repeat this argument starting with $C_\lambda^k(M,g^{\mu,\nu})$ to obtain an inclusion in the other direction, so that in total
$$ C^k_{\lambda}(M, K^{-1} g^{\mu,\nu}) \subset C^k_{\lambda}(M, g'^{\mu,\nu}) \subset C^k_{\lambda}(M, K g^{\mu,\nu}). $$
Taking suprema over the $\cal NG$-dimensions of closed linear subspaces with $K g^{\mu,\nu} = g^{K^{\half} \mu, K^{\half} \nu}$,
\begin{align*}
     \cal G_k(M\to B, \nabla, g)({K^{-\half}}\mu,{K^{-\half}}\nu) 
     &\leq \cal G_k(M\to B, g', \nabla)(\mu,\nu)
     \\&\leq \cal G_k(M\to B, \nabla, g)(K^{\half} \mu,K^{\half}\nu) 
\end{align*}
and hence the spectral density functions are dilatationally equivalent, 
\[G_k(M\to B, \nabla, g)\sim G_k(M\to B, g', \nabla). \qedhere \]
\end{proof}

%% file: FibreHomotopyInvar.tex
Next, we want to study the behaviour of the two-parameter Novikov-Shubin numbers under fibre homotopy equivalences.
Such a homotopy equivalence $f$, say between $M\to B$ and $M'\to B$, should respect the decomposition of $TM \cong_\nabla HM \oplus VM$ and $TM' \cong_{\nabla'} HM' \oplus VM'$ coming from the connections in the sense that $f^*\nabla' = \nabla$. This leads us to the following definition of geometric fibre homotopy equivalences.

\definition{
Let $F_\bullet\to M\xrightarrow{\pi} B$ and $F'_\bullet\to M'\xrightarrow{\pi'} B$ be two fibre bundles over $B$ equipped with connections 
$$ TM\cong_\nabla VM\oplus HM, \qquad TM' \cong_{\nabla'} VM' \oplus HM'.$$
A fibre homotopy equivalence $f\colon M\to M'$ is a homotopy equivalence $f\colon M\to M'$ such that $f$ is a fibre map over the identity $\id_B$ of $B$, that is the diagram 
$$\begin{tikzcd}[ampersand replacement=\&, row sep = 30, column sep = 30]
 M \ar[r, "\pi"] \ar[d, "f"] \& B \ar[d, equal] \\
 M' \ar[r, "\pi'"] \& B,
\end{tikzcd}$$ 
commutes, and so is a homotopy equivalence inverse $g$ of $f$ as well as the homotopy $\Phi\colon M\times [0,1] \to M$ between $gf$ and $\id_M$ at every time $t\in [0,1]$.
We call such a fibre homotopy equivalence $f\colon M\to M'$ geometric if it satisfies $f^*\nabla' = \nabla$. 
}

The property of being geometric implies that the pullback $f^*$ commutes not only with the de Rham differential $d$ itself but also with each of the individual summands we identified earlier.

\lemma{
If $f\colon M\to M'$ is a geometric fibre homotopy equivalence then $f^*$ commutes with the differential $d$ and each of its three summands $d = d^{0,1}+d^{1,0}+d^{2,-1}$.
}
\begin{proof}
Since $f$ is geometric, the fibre homotopy equivalence $f$ restricts fibrewise to homotopy equivalences 
$$f|_{F_b}\colon F_b\xrightarrow{\ \simeq\ } F'_b. $$
and the push-forward $f_*\colon TM \to TM'$ restricts to maps 
$$f_*\colon HM\to HM' \quad \text{and}  \quad f_*\colon VM\to VM'.$$
Therefore, the induced chain homotopy $f^*\colon \Omega^\bullet M'\to \Omega^\bullet M$ restricts under the direct sum decompositions to maps on each $(p,q)$-summand, that is, 
$$ f^*_{p,q}\colon \Omega^p(B,\{\Omega^q(F'_\bullet)\}) \to \Omega^p(B,\{\Omega^q(F_\bullet)\}) $$ 
given on $\alpha_{p,q}\in \Omega^p(B,\{\Omega^q(F'_\bullet)\})$ with $p+q=k$ by 
\begin{align*}
     (f^*_{p,q}\alpha)_{p,q}&(U_1,\dots, U_p)(Z_{p+1}, \dots, Z_{k}) \\
&= (f|_{F_\bullet})^*\left(\alpha_{p,q}(U_1,\dots, U_p)\right)(Z_{p+1}, \dots, Z_{k}) \\
&= \alpha_{p,q}(U_1,\dots, U_p)(df(Z_{p+1}), \dots, df(Z_{k})).
\end{align*}

Recall that the differential $d$ on $ \Omega^p(B,\{\Omega^q(F_\bullet)\})$ splits into the following three summands:
\begin{align*}
(d^{0,1} \alpha)_{p,q+1}&(U_1,\dots, U_p)(Z_{p+1}, \dots, Z_{k+1})\\
&=
\sum_{p+1\leq i\leq k+1} (-1)^{i+1} Z_i(\alpha(U_1,\dots, U_p))(Z_{p+1},\overset{\widehat{Z_i}}{\dots}, Z_{k+1})
\\
&\quad + \sum_{p+1\leq i<j \leq k+1} (-1)^{i+j+1-p} \alpha(U_1,\dots, U_p)([Z_i, Z_j], Z_{p+1}, \overset{\widehat{Z_i},\widehat{Z_j}}{\dots}, Z_{k+1}),
\\
(d^{1,0} \alpha)_{p+1,q}&(U_1,\dots, U_{p+1})(Z_{p+2}, \dots, Z_{k+1})\\
&= 
\sum_{1\leq i\leq p+1} (-1)^{i+1} {\widetilde U_i}(\alpha({U_1},\overset{\widehat{U_i}}{\dots},{U_{p+1}})(Z_{p+2},\dots,Z_{k+1})) \\
&\quad + \sum_{1\leq i<j\leq p+1} (-1)^{i+j+1} \alpha([{U_i},{U_j}],{U_1},\overset{\widehat{U_i},\widehat{U_j}}{\dots},{U_{p+1}})( Z_{p+2},\dots,Z_{k+1}) \\
&\quad + \sum_{1\leq i \leq p+1 < j\leq k+1} (-1)^{i+j+1-p} \alpha({U_1},\overset{\widehat{U_i}}{\dots}, {U_{p+1}})( [\widetilde{U_i},Z_j], Z_{p+2},\overset{\widehat{Z_j}}{\dots},Z_{k+1})
\\
(d^{2,-1} \alpha)_{p+2,q-1}&(U_1,\dots, U_{p+2})(Z_{p+3}, \dots, Z_{k+1})\\
&= \sum_{1\leq i<j\leq p+2} (-1)^{i+j+1-p} \alpha({U_1},\overset{\widehat{U_i},\widehat{U_j}}{\dots}, {U_{p+2}})([\widetilde{U_i},\widetilde{U_j}]_V, Z_{p+3},\dots,Z_{k+1}).
\end{align*}
Here, $f^*$ commutes with $d^{0,1}$ as we can see directly from the formulae or from the fact that 
$$(d^{0,1} \alpha)(U_1,\dots,U_p) = d_{F_\bullet}(\alpha(U_1,\dots,U_p))$$ 
acts as the fibre differential and therefore commutes with the pullback of $f$.
From the formulae we see further that $d^{1,0}$ commutes with $f^*$ since $df(\widetilde{U}) = \widetilde{U}'\circ f = \widetilde{U}$ as $f$ preserves base points, 
$$ df([\widetilde{U}, Z]) = [df(\widetilde{U}), df(Z)] = [\widetilde U', df(Z)], $$ 
where $\widetilde{U}$ is the horizontal lift of $U$ to $TM$ and $\widetilde{U}'$ the horizontal lift to $TM'$.
Lastly, $d^{2,-1}$ commutes with $f^*$ since
\begin{align*}
     df\left([\widetilde{U_1},\widetilde{U_2} ]_V\right) 
     &= df\left([\widetilde{U_1},\widetilde{U_2} ] - [\widetilde{U_1},\widetilde{U_2} ]_H \right) 
     = [df(\widetilde{U_1}),df(\widetilde{U_2}) ] - df(\widetilde{[{U_1},{U_2} ]}) 
     \\
     &= [\widetilde{U_1}', \widetilde{U_2}'] - \widetilde{[U_1,U_2]}'
     = [\widetilde{U_1}', \widetilde{U_2}'] - [\widetilde{U_1}', \widetilde{U_2}']_H
     = [\widetilde U_1', \widetilde U_2']_V. \qedhere
\end{align*}
\end{proof}

\lemma{
Let $G\acts (M\to B, \nabla, g)$ and $G\acts (M'\to B, \nabla', g')$ be two Riemannian fibre bundles with connection over the same base $B$ and with compatible $G$-action.
Let $f\colon M\to M'$ be a $G$-equivariant geometric fibre homotopy equivalence. 
If $f^*$ and a geometric fibre homotopy inverse $g^*$ of $f^*$ are bounded as operators between $L^2\Omega^\bullet M'$ and $L^2\Omega^\bullet M$, then the two-parameter spectral density functions are dilatationally equivalent, that is, for $0\leq k\leq \dim(M)$,
$$\cal G_k(M'\to B, \nabla') \sim \cal G_k(M\to B, f^*\nabla').$$
}
\begin{proof}
By assumption, the induced map $f^*$ is a bounded chain homotopy equivalence $L^2\Omega^\bullet M'\to L^2\Omega^\bullet M$ of Hilbert chain complexes, with bounded inverse $g^*$. 
Since 
$$\cal G_k(M'\to B, \nabla',g')(\mu,\nu) = \cal F_k(M',g'^{\mu,\nu})(1)$$
and in the same way
$$\cal G_k(M\to B, f^*\nabla', g)(\mu,\nu) = \cal F_k(M,g^{\mu,\nu})(1)$$ for some\footnote{By Theorem~\ref{Thm_2NSI_Metric_fixed_conn_invar}, the dilatational equivalence classes of these spectral density functions are independent of this choice.} $G$-invariant Riemannian metrics compatible with the connections, the statement follows from a Proposition of M.~Gromov and M.~A.~Shubin~\cite[Prop.~4.1]{GroShu2}: 

There exists $C(\mu,\nu)$ depending only on $\|f^*\|_{(M',g'^{\mu,\nu})\to (M,g^{\mu,\nu})}$ and $\|g^*\|_{(M',g'^{\mu,\nu})\to (M,g^{\mu,\nu})}$ with
\begin{align*}
    \cal G_k(M\to B, f^*\nabla', g)&(C(\mu,\nu)^{-1}\mu, C(\mu,\nu)^{-1}\nu)
    \\
    &= \cal F_k(M,C(\mu,\nu)^{-1} g^{\mu,\nu})(1)  \\
    &= \cal F_k(M,g^{\mu,\nu})(C(\mu,\nu)^{-1})  
    \\ 
    &\leq \cal F_k(M',g'^{\mu,\nu})(1) \\
    &=\cal G(M'\to B, \nabla', g')(\mu, \nu) \\
    &\leq \cal F_k(M,g^{\mu,\nu})(C(\mu,\nu)) \\
    &= \cal G_k(M\to B, f^*\nabla', g)(C(\mu,\nu)\mu, C(\mu,\nu)\nu).
\end{align*} 
Since for $f^*\colon \Omega^p(B,\{\Omega^q F'_\bullet\})\to \Omega^p(B,\{\Omega^q F_\bullet\})$ (and in the same way for $g^*$), 
\begin{align*}
     \|f^*\|_{(M',g'^{\mu,\nu})\to (M,g^{\mu,\nu})} 
    &= \sup_{0\neq\omega\in \Omega^p(B,\{\Omega^q F'_\bullet\}) } \frac{\|f^*\omega\|_{g^{\mu,\nu}}}{\|\omega\|_{g'^{\mu,\nu}}} 
    \\
    &= \sup_{0\neq\omega\in \Omega^p(B,\{\Omega^q F'_\bullet\}) } \frac{\mu^{-p} \nu^{-q} \cdot \|f^*\omega\|_{g}}{\mu^{-p} \nu^{-q} \cdot \|\omega\|_{g'}}
    \\
    &= \sup_{0\neq\omega\in \Omega^p(B,\{\Omega^q F'_\bullet\}) } \frac{\|f^*\omega\|_{g}}{\|\omega\|_{g'}} 
    = \|f^*\|_{(M',g')\to (M,g)},
\end{align*}
the norms of $f^*$ and $g^*$ are independent of $\mu,\nu$ and hence so is $C=C(\mu,\nu)$. Therefore, the claim follows from the inequalities above.
\end{proof}

Following the idea behind M.~Gromov and M.~A.~Shubin's approach in~\cite{GroShu} further, 
we can drop the restrictive requirement that $f^*$ and $g^*$ are bounded and obtain the desired first invariance theorem.

\theorem{
\label{Thm_2NSI_fibre_homot_invar}
In the notation above, if there is a $G$-equivariant geometric fibre homotopy equivalence between $M\to B$ and $M'\to B$, then for $0\leq k\leq \dim(M)$,
$$\cal G_k(M'\to B, \nabla') \sim \cal G_k(M\to B, f^*\nabla').$$
}

\begin{proof}
In the spirit of \cite[Thm~5.2]{GroShu2}, we show that for any geometric fibre homotopy equivalence $f\colon M\to M'$, we can construct a homotopy equivalence between the corresponding Hilbert chain complexes (which in particular is bounded). The main step here is to construct a submersive fibre homotopy equivalence $\widetilde f\colon M\times D^N\to M'$ from the a thickened fibre bundle $F_\bullet\times D^N\to M\times D^N \to B$ to $F'_\bullet\to M'\to B$, where $D^N$ is a disk in $\R^N$. 

We consider the vertical bundle $VM'\to M'$ and its pullback $f^*VM'\to M$ along $f$. 
By the smooth Serre-Swan theorem\footnote{The original Serre-Swan theorem~\cite[Lem.~5]{Swan} holds for compact topological manifolds. It since has been shown that in the smooth case it holds also for non-compact manifolds, see for example J.~Nestruev's book \cite[Sec.~12.33]{Nestruev} or Section 11.33 in the first edition. 
Here, we want the fibre bundle to be compatible with the action of $G'$ on $B$, so that we may use the Serre-Swan theorem over the compact quotient $f^*V(\leftquotient{G}{M'})\to \leftquotient{G}{M}$ and lift the bundle $\leftquotient{G}{M} \times \R^N \to \leftquotient{G}{M}$ to a bundle $M\times \R^N\to M$ compatible with the group action.
This is possible since $f$ is $G$-equivariant.}
there exists $N\in \N$ and an epimorphism $p_1$ of bundles over $M$, 
$$\begin{tikzcd}[ampersand replacement = \&]
M\times \R^N\ar[r,"p_1"] \ar[d] \& f^*VM' \ar[d] \\
M \ar[r, equal] \& M.
\end{tikzcd}$$
This gives us the following commutative diagram, where $p_1$ and $p_2$ are bundle projections:
$$\begin{tikzcd}[ampersand replacement = \&]
M\times \R^N \ar[rr, twoheadrightarrow, "p_1"]\ar[dr] \& \& f^*VM' \ar[dl]\ar[rr, twoheadrightarrow, "p_2"] \& \& VM' \ar[d, "\pi"] \\
\& M \ar[rrr, "f"]\ar[d] \& \& \& M' \ar[d] \\
\& B \ar[rrr, equal] \& \& \& B
\end{tikzcd}$$
After fixing any $\nabla'$-compatible Riemannian metric $g'$ and the corresponding fibrewise geodesic flows on $M'$, on each fibre $F'_b=\pi^{-1}(b)$ of the bundle $VM'\to M'$ the exponential maps 
$ \exp_b \colon V_bM'\to F_b' $
are defined and they glue to a map 
$$ \exp_V \colon VM'\to F'_\bullet. $$
For each $b\in B$, there is $\varepsilon(b)>0$ such that the exponential map restricts to a diffeomorphism from $D_{\varepsilon(b)}^{V_bM'} = \set{v\in V_bM'}{g'_{V,b}(v,v)< \varepsilon(b)^2}$ onto its image. 
This radius $\varepsilon(b)$ can be chosen to depend continuously on $b$ and be invariant under the cocompact action $G'\acts B$ and such that 
$$ \varepsilon = \inf_{b\in B}\left\{ \varepsilon(b) \right\} = \min_{[b]\in {G'}\backslash{B}} \left\{\varepsilon([b])\right\} $$
exists and $\varepsilon>0$. 
Since $g_{F,b}$ depends smoothly on $b\in B$, the set 
$$ U = \bigcup_{b\in B} D_\varepsilon^{V_bM'} $$
defines a neighbourhood of the zero section $0\in \Gamma(VM')$.
In particular, the map $\exp_V$ restricts to a diffeomorphism from $U$ onto its image in $M'$,
$$ \exp_V \colon U \xrightarrow{\ \cong\ } \exp_V(U). $$
Further, for each $b\in B$ we can find $\delta(b)>0$ depending continuously on $b$ such that the subset $\{b\}\times D_\delta^N$ of the fibre over $b$ of $M\times \R^N\to M$ maps into $D_\varepsilon^{V_bM'}$ via the composition $p_2\circ p_1$. 
Since $f$ preserves the base point, this can be chosen invariantly under the cocompact action $G'\acts B$ and we can define $\delta = \min_{[b]\in G'\backslash B} \left\{\delta([b])\right\} > 0$. 
The image of $M\times D_\delta^N\subset M\times \R^N$ under $p_2\circ p_1$ is contained in $U$.
Hence, the composition 
$$\widetilde{f}=\exp_\bullet\circ \, p_2\circ p_1$$ 
defines a submersion from $M\times D_\delta^N$ into $M'$ (as a map over $\id_B$):
$$
\begin{tikzcd}[ampersand replacement=\&, row sep = 10] 
M\times \R^N \ar[d, phantom, "\cup" description] \& \& VM' \ar[d, phantom, "\cup" description] \&
\\
M\times D_\delta^N \ar[rrrr, "\widetilde f" description, rounded corners, to path ={[pos=0.25]  
  (\tikztostart.south)
    -| ([yshift=-0.5cm]\tikztostart.south)
    -| (\tikztotarget)\tikztonodes
    -- (\tikztotarget.south)}] 
    \ar[r, twoheadrightarrow, "p_1"] \& f^*VM' \ar[r, twoheadrightarrow, "p_2"] \& U  \ar[r, "\exp_V", "\cong"'] \& \exp_V(U)  \ar[r, phantom, "\subset" description] \& M'
\end{tikzcd}
$$
Denote by $\iota\colon M\cong M\times \{0\} \hookrightarrow M\times D^N_\delta$ the inclusion as the zero section. 
Then the following diagram commutes:
$$\begin{tikzcd}[ampersand replacement = \&]
M \ar[r, "\iota"]\ar[d]\ar[rr, "f" description, rounded corners, to path ={[pos=0.25]  
  (\tikztostart.north)
    -| ([yshift=.5cm]\tikztostart.north)
    -| (\tikztotarget.north)\tikztonodes
    -- (\tikztotarget)}] \& M\times D_\delta^N \ar[r, "\widetilde f"]\ar[d] \& M'\ar[d] \\
B \ar[r, equal] \& B \ar[r, equal] \& B 
\end{tikzcd}$$
Note that all maps are bundle maps over the identity $\id_B$.
The cochain homotopy equivalences $L^2 \Omega^k(M) \simeq L^2 \Omega^k(M\times D_\delta^N)$ respect the direct sum decompositions.\footnote{
These homotopy equivalences are explicitly constructed in~\cite[Lem. 5.1]{GroShu2}: Let $I=[0,1]$ and $p\colon M\times I \to M$ be the natural projection and let $i_t\colon M\to M\times  I$ for $t\in I$ be that map $x\mapsto (x,t)$. 
Then $p^*\colon L^2\Omega^k M \to L^2 \Omega^k (M\times I)$ is a homotopy equivalence with inverse $J\colon L^2\Omega(M\times I) \to L^2\Omega^k M$, $J\omega = \int_0^1 i_t^*\omega dt$. Using this and the fact that $I^N\simeq D_\delta^N$ by Lipschitz maps gives the needed homotopy equivalences. Here, we consider $M\times I$ as a bundle $M\times I\to B$ with fibres $F_b\times I$ over $b\in B$.}
Following \cite[Thm.~5.2]{GroShu2} further, the cochain homotopy equivalence $\widetilde f^*$ between $L^2 \Omega^\bullet M'$ and $L^2 \Omega^\bullet (M\times D^N_\delta)$ induced by the submersion $\widetilde f$ is bounded.
Since $f$ is a bundle map over $\id_B$, we even obtain bounded homotopy equivalences on each summand of the direct sum decomposition, $\Omega^p(B, \{\Omega^q F'_\bullet\})\to \Omega^p(B, \{\Omega^q F_\bullet\})$. 
The claim now follows from the previous lemma.
\end{proof}

%% file: MetricInvar.tex
We have seen that the two-parameter Novikov-Shubin numbers behave well if the connection is fixed.
If we allow the connection to vary, we still obtain invariance properties if we scale the fibre at least as fast as the base, that is on the parameter space $\{\nu \leq \mu\}$.

\theorem{
\label{Thm_2NSI_half_metric_indep}
Let $G$ be a group and $M\to B$ be equipped with two pairs of connection and Riemannian metric such that $G\acts (M\to B,\nabla,g)$ and $G\acts (M\to B,\nabla',g')$ are Riemannian fibre bundles with connection and compatible $G$-action. 
Then for all $0\leq k\leq \dim(M)$ the two-parameter spectral density functions restricted to the subspace $\{\nu\leq\mu\}$ are dilatationally equivalent,
$$ \cal G_k(M,\nabla,g)|_{\{\nu\leq\mu\}}\sim \cal G_k(M,\nabla',g')|_{\{\nu\leq\mu\}}. $$ 
}
\begin{proof}
Consider the decompositions 
$$ VM\oplus HM \quad \cong_\nabla \quad TM \quad  \cong_{\nabla'} \quad  VM \oplus H'M, $$
where the vertical bundle $VM = \ker (\pi^*)$ is independent of the connection. 
The identity 
$$\id\colon (M,g)\to (M,g')$$ 
induces a map $d\id \colon TM\to TM$ decomposing into maps $d\id \colon VM\to VM$ and $d\id\colon HM\to VM\oplus H'M$, so vertical tangent vectors remain vertical, but horizontal tangent vectors can obtain a vertical component. 
This is captured in the following diagram.

$$\begin{tikzcd}[ampersand replacement = \&, row sep = 30, column sep = 40]
HM \ar[d, phantom, "\oplus" description] \ar[r, "d\id"' {yshift=-5pt, xshift = 0pt}] \ar[dr] \& H'M \ar[d, phantom, "\oplus" description]
\\
VM \ar[d, phantom, "\cong_\nabla" description] \ar[r, "d\id"] \& VM  \ar[d, phantom, "\cong_{\nabla'}" description]
\\
TM \ar[d] \ar[r, "d\id"] \& TM \ar[d]
\\
(M,g) \ar[r, "\id"] \& (M,g') 
\end{tikzcd}$$

For a form of pure degree $(p,q)$ with respect to the direct sum decomposition coming from the connection $\nabla'$,
$$\omega_{p,q} \in \Omega^p(B,\{\Omega^q(F_\bullet)\}) \subset_{\nabla'} \Omega^k (M,g'), $$
its pullback
$ \id^* \omega_{p,q} \in \Omega^k (M,g) $
decomposes under the direct sum decomposition coming from the connection $\nabla$ as a sum,
$$ \id^* \omega_{p,q}  = \sum_{r+s=k} \alpha_{r,s}, $$ 
with $\alpha_{r,s}\in \Omega^r(B,\{\Omega^s(F_\bullet)\}) \subset_{\nabla} \Omega^k (M,g)$.
Since 
$$\id^* \omega(X_1,\dots, X_k) = \omega(d\id(X_1),\dots, d\id(X_k)),$$
in the $\nabla$-decomposition the $(r,s)$-summand $\alpha_{r,s}$ vanishes if ${r<p}$ or equivalently ${s>q}$. 
Hence
$$ \id^* \omega_{p,q}  = \sum_{\substack{r+s=k \\ r\geq p\ \wedge\  s\leq q}} \alpha_{r,s}. $$
Therefore,
$ \|\omega_{p,q}\|_{g'^{\mu,\nu}} = \mu^{-p} \nu^{-q}\|\omega_{p,q}\|_{g'}$
and
\begin{align*}
     \|\id^* \omega_{p,q}\|_{g^{\mu,\nu}} 
     &= \sum_{\substack{r+s=k \\ r\geq p\ \wedge\ s\leq q}} \|\alpha_{r,s}\|_{g^{\mu,\nu}}
     = \sum_{\substack{r+s=k \\ r\geq p\ \wedge\ s\leq q}} \mu^{-r} \nu^{-s} \|\alpha_{r,s}\|_{g} \\
     &\!\!\overset{\nu\leq\mu} \leq \sum_{\substack{r+s=k \\ r\geq p\ \wedge\ s\leq q}} \mu^{-p} \nu^{-q} \|\alpha_{r,s}\|_{g} 
     = \mu^{-p} \nu^{-q} \|\id^*\omega_{p,q}\|_g.
\end{align*}
Consequently, 
\begin{align*}
    \|\id^*|_{\Omega^p(B,\{\Omega^q(F_\bullet)\})}\|_{(M,g'^{\mu,\nu})\to (M,g^{\mu,\nu})}
    &= 
    \sup_{\omega_{p,q}\in {\Omega^p(B,\{\Omega^q(F_\bullet)\})}} \frac{\|\id^*\omega_{p,q}\|_{g^{\mu,\nu}}}{\|\omega_{p,q}\|_{g'^{\mu,\nu}}} 
    \\
    &\leq
    \sup_{\omega_{p,q}\in {\Omega^p(B,\{\Omega^q(F_\bullet)\})}} \frac{\mu^{-p} \nu^{-q}\cdot \|\id^*\omega_{p,q}\|_{g}}{\mu^{-p} \nu^{-q}\cdot \|\omega_{p,q}\|_{g'}} 
    \\
    &=
    \|\id^*|_{\Omega^p(B,\{\Omega^q(F_\bullet)\})}\|_{(M,g')\to (M,g)}
\end{align*} 
Since the decomposition into the $\Omega^p(B,\{\Omega^q(F_\bullet)\})$ is orthogonal, it follows that 
$$ \|\id^*\|_{(M,g'^{\mu,\nu})\to (M,g^{\mu,\nu})} \leq \|\id^*\|_{(M,g')\to (M,g)}. $$
The same argument holds if we consider the identity map as a map in the other direction, that is $\id\colon (M,g')\to (M,g)$. 
Let 
$$K=\max\left\{\|\id^*\|_{(M,g')\to (M,g)},\, \|\id^*\|_{(M,g)\to (M,g')}\right\}.$$
For any $\omega\in C^k(M,g'^{\mu,\nu})(1)$ with $\nu\leq\mu$ it follows, therefore, that
\begin{align*}
\| d\id^*\omega \|_{g^{\mu,\nu}}
&= \| \id^* d\omega\|_{g^{\mu,\nu}} 
\leq K \|d\omega\|_{g'^{\mu,\nu}}
\leq K \|\omega\|_{g'^{\mu,\nu}}
= K \|\id^* \omega\|_{g'^{\mu,\nu}} 
\leq K^2 \|\omega\|_{g^{\mu,\nu}} 
\end{align*}
and similarly in the other direction.
These inequalities imply that 
$$ C^k(M,g^{\mu,\nu})(K^{-2}) \subset C^k(M,g'^{\mu,\nu})(1) \subset  C^k(M,g^{\mu,\nu})(K^2). $$
Hence the spectral density functions are dilatationally equivalent and the claim follows.
\end{proof}

%% file: Example.tex
Let us consider the three-dimensional Heisenberg group 
$\mathbb H^3$ 
with its associated Lie algebra 
$\frak h^3 = \langle X,Y,Z \:|\: [X,Y]=Z\rangle$
as a fibre bundle
with fibre $\R$ corresponding to the central $Z$-direction and base $\R^2$ corresponding to the $X$- and $Y$-directions.
A basis of left-invariant vector fields is given by the vector fields
\begin{align*}
\vartheta_X &= \partial_X - \frac{1}{2} y \partial_Z, \qquad
\vartheta_Y = \partial_Y + \frac{1}{2} x \partial_Z, \qquad
\vartheta_Z = \partial_Z,
\end{align*}
where $x$ and $y$ denote coordinates in the base $\R^2 = \langle X,Y\rangle$.

Requiring that $\vartheta_X, \vartheta_Y$ and $\vartheta_Z$ are orthonormal yields the standard metric $g$ and with 
$VM = \langle \vartheta_Z\rangle$ and $HM = \langle \vartheta_X,\vartheta_Y\rangle$. 
We also fix a connection $\nabla$. 
The scaled metric $g^{\overline\mu,\overline\nu}$ is the metric for which 
$${\overline\mu}^{-1}\cdot\vartheta_X, \quad {\overline\mu}^{-1}\cdot\vartheta_Y \quad \text{ and } \quad {\overline\nu}^{-1}\cdot\vartheta_Z$$ 
form an orthonormal basis of $\frak h^3$. 
Refining a computation of J.~Lott~\cite[Prop.~52]{Lott}, we obtain the following values for the two-parameter Novikov-Shubin numbers.

\begin{figure}[h!]
\centering
\begin{minipage}{.5\textwidth}
\begin{tikzpicture}[scale=0.85]
    \begin{axis}[xmin=-1, xmax=1, ymin=2.75, ymax=6.25, xlabel={$\zeta$}, ylabel={$\alpha_0(\frak h^3)(\lambda,\lambda^{1+\zeta})$}]
        \addplot[domain = -0.5:1] {4+2*x};
        \addplot[domain = -1:0, dashed] {4};
        \addplot [dashed] coordinates {(0, 2.75) (0, 4)};
        \addplot[mark=*] coordinates {(0,4)};
    \end{axis}
\end{tikzpicture}
\end{minipage}\hfill%
\begin{minipage}{.5\textwidth}
\begin{tikzpicture}[scale=0.85]
    \begin{axis}[xmin=-1, xmax=1, ymin=-0.25, ymax=3.25, xlabel={$\zeta$}, ylabel={$\alpha_1(\frak h^3)(\lambda,\lambda^{1+\zeta})$}]
        \addplot[domain = -0.5:1] {2-2*x};
        \addplot[domain = -1:0, dashed] {2};
        \addplot [dashed] coordinates {(0,-0.25) (0, 2)};
        \addplot[mark=*] coordinates {(0,2)};
    \end{axis}
\end{tikzpicture} 
\end{minipage}%
\hfill
\caption{The two-parameter Novikov-Shubin numbers of $\frak h^3$. On the left, we see a plot for $\alpha_0(\frak h^3)(\lambda,\lambda^{1+\zeta})$ and on the right for $\alpha_1(\frak h^3)(\lambda,\lambda^{1+\zeta})$. 
The marked points at $\zeta=0$ indicate the classical Novikov-Shubin invariants $\alpha_0(\mathbb H^3)$ and $\alpha_1(\mathbb H^3)$. 
For $\alpha_0$, the contributions of base and fibre seem to agree. 
In particular, as $\zeta$ increases, so does $\alpha_0$. 
For $\alpha_1$, the opposite is the case: As $\zeta$ increases, $\alpha_1$ decreases. 
This gives an interesting insight to (classical) Novikov-Shubin invariants. 
Comparing the Novikov-Shubin invariants for $\mathbb H^3$ and $\R^3$, 
$\alpha_0(\mathbb H^3) = 4 > 3 = \alpha_0(\R^3)$ but $\alpha_1(\mathbb H^3) = 2 < 3 = \alpha_1(\R^3)$.
This fits to the observation that in $\mathbb H^3$, the $Z$-direction scales like the product of the base directions ($[aX, bY] = ab Z$), so scaling the fibre with $\lambda^{\half}$ seems to counteract this.
}
\label{F_2PNSIForH3}
\end{figure}
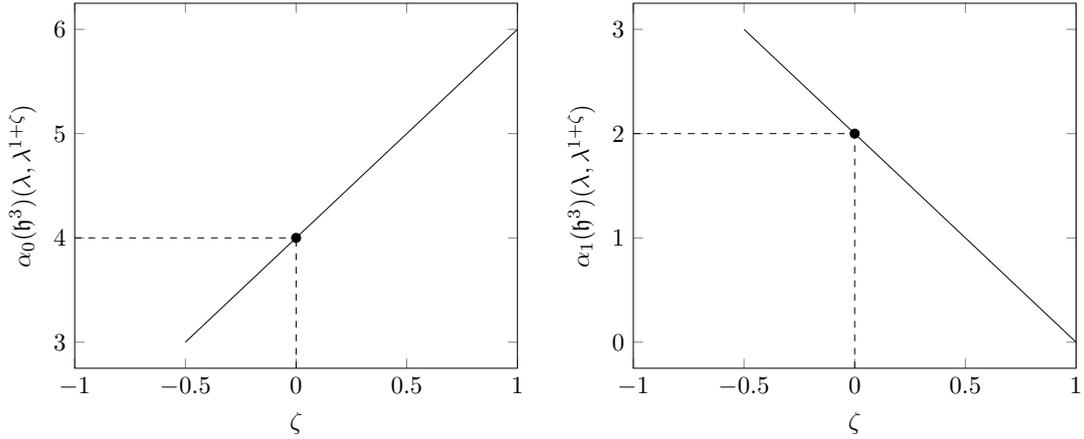

\theorem{
\label{T_result2ParamH3}
On $\frak h^3$, by direct computation we obtain
\begin{align*}
    \alpha_0(\frak h^3)(\lambda,\lambda^{1+\zeta}) &= 4+2\zeta \quad \text{for } -\half \leq \zeta, \\
     \alpha_1(\frak h^3)(\lambda,\lambda^{1+\zeta}) &=  2-2\zeta \quad \text{for } -\half <  \zeta < 1,
\end{align*}
and, by Hodge duality, also 
$\alpha_2(\frak h^3)(\lambda,\lambda^{1+\zeta}) = 4+2\zeta$ for $-\half \leq \zeta$. Compare also Figure~\ref{F_2PNSIForH3}. 
}

\begin{proof}
It was shown by J. Lott~\cite[Prop.~52]{Lott} that in this setting of $\mathbb H^3$ with metric $g^{1,c}$, the heat kernel on functions is given by
\begin{align*}
e^{-t\Delta_0}(0,0) = \frac{1}{4\pi^2}\frac{1}{ct^2} \int_{0}^\infty  e^{-\frac{u^2}{c^2t}} \sinh(u)^{-1} u \,\mathrm{d}u.
\end{align*}
Classically, if $c$ is constant and we let $t\to \infty$, the density function of the normal distribution $e^{-\frac{u^2}{c^2t}}$, converges to the constant-$1$ function and therefore
$$\lim_{t\to\infty} \int_{0}^\infty  e^{-\frac{u^2}{c^2t}} \sinh(u)^{-1} u \,\mathrm{d}u = \int_0^\infty  \sinh(u)^{-1} u \,\mathrm{d}u = \frac{\pi}{4}. $$
Hence, $e^{-t\Delta_0}(0,0)$ is in $\Theta (t^{-2})$ as $t\to \infty$ and, following M.~Gromov and M.~A.~Shubin's work in~\cite{GroShu2}, we can therefore conclude that $\alpha_0(\mathbb H^3)=4$. \smallskip

If we let $c$ depend on $t$, the same argument remains true as long as $c(t)^2 t\to \infty$ as $t\to \infty$, showing that then 
$$e^{-\Delta_0}(0,0)\in \Theta (c(t)^{-1}t^{-2}).$$
Therefore, with $c = t^\zeta$ and $\zeta > -1/2$,
$$ \alpha_0(\frak h^3)(\lambda, \lambda^{1+\zeta}) = \alpha\left(\lambda\mapsto \cal G_{0}(\lambda, \lambda^{1+\zeta})\right) = 4+2\zeta.$$
Indeed, since for $\zeta = -\half$ the integral is a positive constant,
$$ 0 < \int_{0}^\infty  e^{-{u^2}} \sinh(u)^{-1} u \,\mathrm{d}u < \infty, $$
the argument holds also for $\zeta = -\half$, however, the integral converges to zero for $\zeta < -\half$, so that its asymptotic behaviour needs to be taken into account.
The summand $2\zeta$ tells us that the scaling of the $Z$-direction contributes quadratically to the spectral density. 
This fits with the computation of $\alpha_0(\mathbb H^3)=N(\mathbb H^3)$ via the growth rate $N(\mathbb H^3)$ (using the result of N.~Th.~Varopolous~\cite{Varop}) since by the Bass-Guivarc'h formula, 
\begin{align*}
N(\mathbb H^3) &= \rk(\langle X,Y\rangle) + 2\cdot \rk(\langle Z\rangle) 
= 2 + 2 = 4,
\end{align*}
so we also see a quadratic contribution from the central $Z$-direction in this picture. \smallskip

On $1$-forms, J.~Lott computes the heat operator as
\begin{align*}
e^{-t\Delta_1}(0,0) &= \frac{1}{2\pi^2}\frac{1}{c}
\left[ I_1^+ + I_1^- + I_2 + I_3 \right], 
\end{align*}
where the summands $I_\bullet$ are the following integral expressions:
\begin{align*}
    I_1^\pm &= \int_0^\infty \sum\limits_{m=1}^\infty e^{-t\left[
(2m+1)k+ \frac{k^2}{c^2}+\frac{c^2}{2}   \pm c\sqrt{(2m+1)k+\frac{k^2}{c^2}+\frac{c^2}{4}}\right]} kdk,
\\
I_2 &= \int_0^\infty e^{-\frac{k^2}{c^2}t} kdk, 
\\
I_3 &= \int_0^\infty  e^{-\left(2k+\frac{k^2}{c^2}+c^2\right)t} kdk.
\end{align*}
J.~Lott estimates these integrals in the case where $c$ is constant in order to compute the Novikov-Shubin invariant $\alpha_1(\mathbb H^3)=2$.  
We will now adapt these computations to the case where $c=c(t)$ is a function of $t$, in particular a power $c(t)=t^\zeta$.
\lemma{The integrals $I_2$ and $I_3$ evaluate to
\begin{align*}
    I_2 &= \frac{1}{2}\frac{c^2}{t}, \\
    I_3 &= \frac{1}{2}\frac{c^2}{t}e^{-c^2 t} + \sqrt{\pi}\frac{c^3}{\sqrt{t}} \cdot \mathrm{erfc}\left(c\sqrt{t}\right),
\end{align*}
where $\mathrm{erfc}$ denotes the complementary Gauss error function.
}
\begin{proof}
The integral $I_2$ can be directly evaluated by substituting $u = k^2$ as 
\begin{align*}
    I_2 &= \int_0^\infty e^{-\frac{t}{c^2}\cdot k^2} k \,\mathrm{d}k 
    = \frac{1}{2} \int_0^\infty e^{\frac{-t}{c^2} u} \,\mathrm{d}u 
    = \frac{1}{2}\frac{c^2}{t}.
\end{align*}

Substituting $u=(\nicefrac{k}{c}+c)^2t$ and $v=(\nicefrac{k}{c}+c)\sqrt{t}$, we can compute
\begin{align*}
    I_3 &= \int_0^\infty  e^{-\left(2k+\frac{k^2}{c^2}+c^2\right)t} k \,\mathrm{d}k 
    \\
    &= \int_0^\infty  e^{-(\frac{k}{c}+c)^2t} k \,\mathrm{d}k 
    \\
    &= \frac{c^2}{2t}\int_0^\infty e^{-(\frac{k}{c}+c)^2t}  \left(\frac{2t}{c^2}k + 2t\right)\,\mathrm{d}k - c^2 \int_0^\infty e^{-(\frac{k}{c}+c)^2t} \,\mathrm{d}k
    \\
    &= \frac{c^2}{2t}\int_{c^2t}^\infty e^{-u}  \,\mathrm{d}u - \frac{c^3}{\sqrt{t}} \int_{c\sqrt{t}}^\infty e^{-{v^2}} \,\mathrm{d}v
    \\
    &= 
    \frac{c^2}{2t}e^{-c^2 t} + \frac{\sqrt\pi c^3}{\sqrt{t}} \cdot \mathrm{erfc}\left(c\sqrt{t}\right). \qedhere
\end{align*}
\end{proof}

\lemma{
By substitution, 
$$I_1^\pm = c^4  \int_0^\infty \left(v\mp \frac{1}{2}\right) e^{-tc^2v^2} \sum_{m=1}^\infty \left[ 1 - \left(\sqrt{1+\frac{(v\mp \frac{1}{2})^2 - \frac{1}{4}}{(m+\frac{1}{2})^2}}\right)^{-1} \right] \,\mathrm{d}v.$$
}
\begin{proof}
Following J. Lott's computations, we substitute in the same way
\begin{align*}
   u_\pm &= \sqrt{(2m+1)k + \frac{k^2}{c^2} + \frac{c^2}{4}} \pm \frac{c}{2} \\ u_\pm^2 &= (2m+1)k + \frac{k^2}{c^2} + \frac{c^2}{2} \pm c\sqrt{(2m+1)k+\frac{k^2}{c^2} +  \frac{c^2}{4}} \\
    k_\pm &= c\sqrt{u_\pm^2 \mp u_\pm c + c^2(m+\nicefrac{1}{2})^2} - \left(m+\frac{1}{2}\right)c^2 \\
    \frac{\,\mathrm{d}k_\pm}{\,\mathrm{d}u_\pm} &= \frac{c(u_\pm \mp\nicefrac{c}{2})}{\sqrt{u_\pm^2 \mp u_\pm c + c^2(m+\nicefrac{1}{2})^2}}.
\end{align*} 
Omitting the index $\pm$ in notation\footnote{For $I_1^+$, the index $+$ is to be used and for $I_1^-$ the index $-$ is to be used.}, we use this with $v=\nicefrac{u}{c}$ to obtain
\begin{align*}
    I_1^\pm &= \int_0^\infty \sum\limits_{m=1}^\infty e^{-t\left[
(2m+1)k+ \frac{k^2}{c^2}+\frac{c^2}{2}   \pm c\sqrt{(2m+1)k+\frac{k^2}{c^2}+\frac{c^2}{4}}\right]} k \,\mathrm{d}k
\\
&=c^2 \int_0^\infty \left(u\mp \frac{c}{2}\right) e^{-tu^2} \sum_{m=1}^\infty \frac{\sqrt{u^2 \mp cu + c^2(m+\frac{1}{2})^2} - c(m+\frac{1}{2})}{\sqrt{u^2 \mp cu + c^2(m+\frac{1}{2})^2}} \,\mathrm{d}u
    \\
    &= c^3\int_0^\infty \left(\frac{u}{c}\mp \frac{1}{2}\right) e^{-t\frac{u^2}{c^2}c^2} \sum_{m=1}^\infty \frac{\sqrt{\frac{u^2}{c^2} \mp \frac{u}{c} + (m+\frac{1}{2})^2} - (m+\frac{1}{2})}{\sqrt{\frac{u^2}{c^2} \mp \frac{u}{c} + (m+\frac{1}{2})^2}} \,\mathrm{d}u
    \\
    &= c^4  \int_0^\infty \left(v\mp \frac{1}{2}\right) e^{-tc^2v^2} \sum_{m=1}^\infty \frac{\sqrt{v^2 \mp v + (m+\frac{1}{2})^2} - (m+\frac{1}{2})}{\sqrt{v^2 \mp v + (m+\frac{1}{2})^2}} \,\mathrm{d}v
    \\
    &= c^4  \int_0^\infty \left(v\mp \frac{1}{2}\right) e^{-tc^2v^2} \sum_{m=1}^\infty \left[ 1 - \left(\sqrt{1+\frac{(v\mp \frac{1}{2})^2 - \frac{1}{4}}{(m+\frac{1}{2})^2}}\right)^{-1} \right] \,\mathrm{d}v. \qedhere
\end{align*}
\end{proof}

\lemma{
We can estimate $I_1^-$ by
\begin{align*}
  \frac{1}{5}\left( \frac{\sqrt \pi}{4} \frac{c}{\sqrt{t^3}} + \frac{1}{4}\frac{c^2}{t}\right) \leq I_1^- \leq \frac{\sqrt \pi}{4} \frac{c}{\sqrt{t^3}} + \frac{1}{4}\frac{c^2}{t}.
\end{align*}
}
\begin{proof}
Consider the function $f\colon \R_{\geq 0}\to \R$ describing the summands,
$$f(x) = 1 - \left(\sqrt{1+\frac{(v\mp \frac{1}{2})^2 - \frac{1}{4}}{(x+\frac{1}{2})^2}}\right)^{-1}.$$
This function is positive, monotonously decreasing with values $f(0) = 1 - (2v+ 1)^{-1}$ and $\lim_{x\to\infty} f(x) = 0$. 
We can therefore estimate the sum over the $f(n)$ by integrals,
\begin{align*}
    \int_2^\infty f(x) \,\mathrm{d}x \leq  \sum_{n=1}^\infty f(n) \leq \int_1^\infty f(x) \,\mathrm{d}x
\end{align*}
To compute these integrals, let $w = (v+\frac{1}{2})^2 - \frac{1}{4}$, 
then
\begin{align*}
    F(x) = \int f(x) dx 
    &= \int 1 - \left(\sqrt{1+\frac{w}{(x+\frac{1}{2})^2}}\right)^{-1} \,\mathrm{d}x
    \\
    &= x - \left(x+\frac{1}{2}\right)\sqrt{1 + \frac{w}{(x+\frac{1}{2})^2}} +  \mathrm{const}
\end{align*}
and we can compute the values 
\begin{align*}
F(1) &= 1-\sqrt{(v+\nicefrac{1}{2})^2 + 2} +  \mathrm{const} ,\qquad F(2) = 2-\sqrt{(v+\nicefrac{1}{2})^2 + 6} +  \mathrm{const}
\end{align*} 
as well as $\lim_{x\to\infty} F(x) = -\nicefrac{1}{2}+  \mathrm{const}$.
Hence, we get bounds on the sum by
\begin{align*}
    \sqrt{\left(v+\frac{1}{2}\right)^2 + 6}-\frac{5}{2} \leq \sum_{m=1}^\infty f(m) \leq \sqrt{\left(v+\frac{1}{2}\right)^2 + 2} - \frac{3}{2}.
\end{align*}
For the lower bound, observe that $g\colon \R_{\geq 0}\to \R$, $v\mapsto \sqrt{\left(v+\nicefrac{1}{2}\right)^2 + 6}-\nicefrac{5}{2}$ satisfies $g(0) =0$, 
$$ g'(v) = \frac{v+\frac{1}{2}}{\sqrt{\left(v+\frac{1}{2}\right)^2+6}}, \quad g''(v) = \frac{6}{((v+\frac{1}{2})^2+3)^{3/2}}, $$
so $g''>0$ meaning that $g'$ is strictly monotonously increasing and has its minimum at $g'(0) = \nicefrac{1}{5}$. This implies $g(v) \geq \nicefrac{v}{5}$.
For the upper bound, we do the same analysis and find that for $h(v) = \sqrt{\left(v+\nicefrac{1}{2}\right)^2 + 2} - \nicefrac{3}{2}$ we have $h(0)=0$ and $h'(v)\leq \lim_{v\to\infty} h'(v) = 1$ implying that $h(v)\leq v$.
Hence we get new bounds 
\begin{align*}
    \frac{v}{5} \leq \sum_{m=1}^\infty f(m) \leq v.
\end{align*}
 Using these bounds, we get bounds on $I_1^-$ by evaluating
\begin{align*}
   c^4  \int_0^\infty \left(v+\frac{1}{2}\right) e^{-tc^2v^2}v  \,\mathrm{d}v
    &= c^4 \int_0^\infty v^2 e^{-tc^2v^2}  \,\mathrm{d}v + \frac{c^4}{2} \int_0^\infty v e^{-tc^2v^2} \,\mathrm{d}v
\end{align*}
By partial integration and with $\kappa = c v \sqrt t $, the first summand is given by
\begin{align*}
    c^4\int_0^\infty v^2 e^{-tc^2v^2}  \,\mathrm{d}v 
    &=  c^2 \left[-\frac{ve^{-tc^2v^2}}{2t}\right]_{v=0}^\infty + \frac{c^2}{2t} \int_0^\infty e^{-tc^2v^2} \,\mathrm{d}v
    \\
    &= 0 + \frac{c}{2\sqrt{t^3}} \int_0^\infty e^{-\kappa^2} \,\mathrm{d}\kappa
    = \frac{\sqrt \pi c}{4\sqrt{t^3}}
\end{align*}
and with $\xi = tc^2v^2$ the second summand is 
\begin{align*}
    \frac{c^4}{2}\int_0^\infty v e^{-tc^2v^2} \,\mathrm{d}v 
    &= \frac{c^2}{4t}\int_0^\infty e^{-\xi} \,\mathrm{d}\xi  = \frac{c^2}{4t}.
\end{align*}
Therefore,
\begin{align*}
  \frac{1}{5}\left( \frac{\sqrt \pi c}{4\sqrt{t^3}} + \frac{c^2}{4t}\right) &\leq I_1^- \leq \frac{\sqrt \pi c}{4\sqrt{t^3}} + \frac{c^2}{4t}. \qedhere
\end{align*}
\end{proof}

\lemma{Let $I_4$ be the part of $I_1^+$ starting at $1$, that is 
$$I_4 = c^4  \int_1^\infty \left(v-  \frac{1}{2}\right) e^{-tc^2v^2} \sum_{m=1}^\infty \left[ 1 - \left(\sqrt{1+\frac{(v- \frac{1}{2})^2 - \frac{1}{4}}{(m+\frac{1}{2})^2}}\right)^{-1} \right] \,\mathrm{d}v. $$
Then 
\begin{align*}
     &\frac{1}{5}\left[-\frac{c^2}{4t}e^{-tc^2}+ \frac{\sqrt\pi}{4}\left(\frac{c}{\sqrt{t^3}} + \frac{c^3}{\sqrt t}\right) \mathrm{erfc}(c\sqrt{t})\right]
     \\
	&\quad \leq I_4 \leq \left[-\frac{c^2}{4t}e^{-tc^2}+ \frac{\sqrt\pi}{4}\left(\frac{c}{\sqrt{t^3}} + \frac{c^3}{\sqrt t}\right) \mathrm{erfc}(c\sqrt{t})\right]. 
\end{align*}
}
\begin{proof}
Similar as for $I^1_-$, in the case of $I_1^+$ we consider 
$$f(x) = 1 - \left(\sqrt{1+\frac{(v - \frac{1}{2})^2 - \frac{1}{4}}{(x+\frac{1}{2})^2}}\right)^{-1}.$$
If $v>1$, the function $f$ is again monotonously decreasing
and we can estimate 
$$ 
   \frac{v-1}{5} \leq  \sqrt{\left(v-\frac{1}{2}\right)^2 + 6}-\frac{5}{2} \leq \sum_{m=1}^\infty f(m) \leq \sqrt{\left(v-\frac{1}{2}\right)^2 + 2} - \frac{3}{2} \leq v-1.
$$
Therefore, we can bound the $(v>1)$-part $I_4$ of $I_1^+$ by evaluating 
\begin{align*}
 \widetilde{I_4} &= c^4  \int_1^\infty \left(v-\frac{1}{2}\right) e^{-tc^2v^2}(v-1)  \,\mathrm{d}v
    \\
    &= c^4 \int_1^\infty \left( v^2-\frac{3}{2}v+\frac{1}{2} \right)e^{-tc^2v^2} \,\mathrm{d}v
    \\
    &=
    c^4 \int_1^\infty v^2 e^{-tc^2v^2} \,\mathrm{d}v 
    - \left(\frac{3}{2}c^4\right) \int_1^\infty v e^{-tc^2v^2} \,\mathrm{d}v 
    +\frac{c^4}{2} \int_1^\infty e^{-tc^2v^2} \,\mathrm{d}v 
    \\
    &= c^2 \left[-\frac{ve^{-tc^2v^2}}{2t}\right]_{v=1}^\infty + \frac{c}{2\sqrt t}
    \int_{c\sqrt{t}}^\infty e^{-\kappa^2} \,\mathrm{d}\kappa
     - \left(\frac{3c^2}{4t}\right) \int_{tc^2}^\infty e^{-\xi} \,\mathrm{d}\xi 
     + \frac{c^3}{2\sqrt{t}} \int_{c\sqrt{t}}^\infty e^{-\kappa^2} \,\mathrm{d}\kappa
    \\
    &= \frac{c^2 e^{-tc^2}}{2t} + \frac{\sqrt \pi c}{4\sqrt{t^3}} \mathrm{erfc}(c\sqrt{t}) - \left(\frac{3c^2}{4t}\right)e^{-tc^2} 
    + \frac{\sqrt{\pi} c^3}{4\sqrt{t}}\mathrm{erfc}(c\sqrt{t})
    \\
    &= -\frac{c^2}{4t}e^{-tc^2}+ \frac{\sqrt\pi}{4}\left(\frac{c}{\sqrt{t^3}} + \frac{c^3}{\sqrt{t}}\right) \mathrm{erfc}(c\sqrt{t}) 
\end{align*}
with $\nicefrac{1}{5}\cdot \widetilde{I_4}  \leq I_4 \leq \widetilde{I_4}$.
\end{proof}

It remains to estimate 
$$ I_5 = c^4  \int_0^1 \left(v-  \frac{1}{2}\right) e^{-tc^2v^2} \sum_{m=1}^\infty \left[ 1 - \left(\sqrt{1+\frac{(v- \frac{1}{2})^2 - \frac{1}{4}}{(m+\frac{1}{2})^2}}\right)^{-1} \right] dv. $$

\lemma{
There is some constant $-\infty<-K<0$ such that
$$ -K \left(\frac{1}{2}\frac{c^2}{t}e^{-tc^2} - \frac{\sqrt{\pi}}{4}\frac{c^3}{\sqrt{t}}\mathrm{erfc}(c\sqrt{t})\right) \leq I_5 \leq 0. $$
}
\begin{proof}
Note that the summands are non-positive and
\begin{align*}
    \left[ 1 - \left(\sqrt{1-\frac{1}{4(m+\frac{1}{2})^2}}\right)^{-1} \right]
    \leq \left[ 1 - \left(\sqrt{1+\frac{(v- \frac{1}{2})^2 - \frac{1}{4}}{(m+\frac{1}{2})^2}}\right)^{-1} \right] \leq 0
\end{align*}
so that  
\begin{align*}
   \sum_{m=1}^\infty \left[ 1 - \left(\sqrt{1-\frac{1}{4(m+\frac{1}{2})^2}}\right)^{-1} \right] \cdot c^4  \int_0^1 \left(v-  \frac{1}{2}\right) e^{-tc^2v^2}  dv \leq  I_5 \leq 0.
\end{align*}
The sum converges to some constant $-\infty < -K < 0$
while 
\begin{align*}
     c^4 \int_0^1 \left(v-  \frac{1}{2}\right) e^{-tc^2v^2}  \,\mathrm{d}v 
      &= \frac{c^2}{2t} \int_0^{tc^2} e^{-\xi} \,\mathrm{d}\xi - \frac{c^3}{2\sqrt{t}}\int_0^{c\sqrt{t}}e^{-\kappa^2} \,\mathrm{d}\kappa \\
      &= \frac{c^2}{2t} - \frac{c^2}{2t}e^{-tc^2} - \frac{\sqrt{\pi} c}{4\sqrt{t}}\mathrm{erfc}(c\sqrt{t}). \qedhere
\end{align*}
\end{proof}

\corollary{
If $c = t^\zeta$ for $\zeta > -\half$, then
$$ e^{-t\Delta_1}(0,0) \sim \frac{c}{t}  = \frac{1}{t^{1-\zeta}} $$
as $t\to\infty$. In particular, for $-\half < \zeta < 1$,
\begin{align*}
    \alpha\left(\lambda \mapsto \cal G_{1} (\lambda, \lambda^{1+\zeta})\right) =  2-2\zeta.
\end{align*}
}
\begin{proof}
The assumption
$\zeta>-\half$ implies $c^2t\xrightarrow{t\to\infty}\infty$ and both $e^{-tc^2}$ and $\mathrm{erfc}(c\sqrt{t})$ decay exponentially. By the previous computations, 
$$ e^{-t\Delta_1}(0,0) \sim \frac{1}{c} \left[I_1^+ + I_1^- + I_2 + I_3\right] \sim \frac{c}{t} + \frac{1}{t^{\nicefrac{3}{2}}}  $$
as $t\to\infty$. The assumption $\zeta > -\half$ implies $t^{-\nicefrac{3}{2}}\in \cal O(\nicefrac{c}{t})$.
In particular, since $\nicefrac{c}{t} = t^{\zeta-1}$, this decays to zero as $t\to\infty$ for $\zeta < 1$.
\end{proof}
This concludes the computation of the asymptotics for $\alpha_\bullet(\frak h^3)(\lambda, \lambda^{1+\zeta})$.
\end{proof}